\documentclass[12pt]{article}
\usepackage{kz}
\usepackage{enumerate}
\usepackage{graphicx}
\usepackage{authblk}
\usepackage[font=scriptsize,labelfont=bf]{caption}
\RequirePackage[colorlinks,linkcolor=blue,citecolor=blue,urlcolor=blue]{hyperref}
\usepackage[plain,noend]{algorithm2e}
\def\@algocf@capt@plain{top}

\begin{document}

\def\spacingset#1{\renewcommand{\baselinestretch}%
{#1}\small\normalsize} \spacingset{1}

%%%%%%%%%%%%%%%%%%%%%%%%%%%%%%%%%%%%%%%%%%%%%%%%%%%%%%%%%%%%%%%%%%%%%%%%%%%%%%

\title{\textbf{The Binary Expansion Randomized Ensemble Test}}
\author[1]{Duyeol Lee}
\author[1]{Kai Zhang} 
\author[1,2]{Michael R. Kosorok}
\affil[1]{Department of Statistics and Operations Research, University of North Carolina, Chapel Hill, NC 27599}
\affil[2]{Department of Biostatistics, University of North Carolina, Chapel Hill, NC 27599}
\maketitle
  
\bigskip
\begin{abstract}
Recently, the binary expansion testing framework was introduced to test the independence of two continuous random variables by utilizing symmetry statistics that are complete sufficient statistics for dependence. We develop a new test based on an ensemble approach that uses the sum of squared symmetry statistics and distance correlation. Simulation studies suggest that this method improves the power while preserving the clear interpretation of the binary expansion testing. We extend this method to tests of independence of random vectors in arbitrary dimension. Through random projections, the proposed binary expansion randomized ensemble test transforms the multivariate independence testing problem into a univariate problem. Simulation studies and data example analyses show that the proposed method provides relatively robust performance compared with existing methods.
\end{abstract}

\noindent%
\vfill

\newpage
\spacingset{1.68} % DON'T change the spacing!

\section{Introduction}

Nonparametric testing of independence is a fundamental problem in statistics and has been studied carefully by many classical papers such as \citet{hoeffding1948non}. This problem has been gaining greater interest recently due to its important roles in machine learning and big data analysis. Some important recent developments include \citet{szekely2007measuring, gretton2008kernel, heller2012consistent, heller2016consistent, heller2016multivariate, han2017distribution, wang2017generalized, pfister2018kernel, gorsky2018multifit, ma2019fisher, deb2019multivariate, berrett2019nonparametric}. \citet{josse2016measuring} have published an authoritative review.

One important problem in nonparametric dependence detection is nonuniform consistency, which means that no test can uniformly detect all forms of dependency, as described by \citet{zhang2019bet}. This problem is particularly severe for nonlinear relationships, which are common in many areas of science. To avoid the power loss due to nonuniform consistency, \citet{zhang2019bet} considers the binary expansion statistics (BEStat) framework; this framework examines dependence with a filtration approach induced by the binary expansion of uniformly distributed variables. \citet{zhang2019bet} also proposed testing independence of two continuous variables with the framework of maximum binary expansion testing (BET). The BET achieves uniform consistency and is minimax optimal in power. It also provides clear interpretability, and it can be implemented efficiently by bitwise operations.

Although the BET works well in testing independence between two variables, two crucial improvements are needed for greater practical applicability. The first requirement is to improve the power of the BET under certain cases such as linear dependency. The second requirement is to extend the test for testing independence of random vectors. In this paper, we describe a new approach that solves both of these problems. The first problem is addressed by a novel ensemble approach, and the second problem is solved by using one-dimensional random projecting. Due to utilization of both a random projection and an ensemble approach, we call the new method the binary expansion randomized ensemble test (BERET). We show with simulation studies that the proposed method has good power properties.

Through example datasets, we illustrate how the proposed method provides clear interpretability while maintaining good power properties across various dependence structures, including both linear and nonlinear relationships. In a life expectancy example, our method is able to detect three meaningful and interpretable relationships and provide similar p-values as competing methods. In a mortality rate example, we show that the canoncial correlation test can be interpretable but fails to detect nonlinear dependence structure. This is unfortunate since the canonical correlation test is the only other method besides the proposed method that has inherent interpretability. In contrast, our method is able to identify meaningful relationships even when there is a nonlinear relationship. In the house price example, the mutual 
 information test fails to reject independence because the linear relationship is not sufficiently strong. However, our method rejects independence because of its boosted sensitivity to linear relationships, and it is also able to detect interpretable dependence structures including linear relationships. The canonical correlation test also works here and also provides interpretability. However, our method is the only method that has power for detecting both linear and nonlinear relationships as well as being able to illuminate interpretable dependency structures.

This article is organized as follows. Section \ref{method}  describes the ensemble method and the BERET procedure. In Section \ref{simulation}, we present simulation studies on performance of proposed method. Section \ref{data} illustrates our method with three data examples. Concluding remarks are presented in Section \ref{conclusion}.

\section{Proposed Method}\label{method}

\subsection{The Binary Expansion Testing Framework}

We briefly introduce the testing procedure and useful notations from \citet{zhang2019bet}. Let $(X_1, Y_1), \dots, (X_n, Y_n)$ be a random sample from distributions of $X$ and $Y$. If the marginal distributions of $X$ and $Y$ are known, we can use the CDF transformation so that $U=F_X(X)$ and $V=F_Y(Y)$ are each uniformly distributed over $[0,1]$. The binary expansions of two random variables $U$ and $V$ can be expressed as $U=\sum_{k=1}^{\infty}{A_k/2^k}$ and $V=\sum_{k=1}^{\infty}{B_k/2^k}$ where $A_k \overset{i.i.d.}{\sim} Bernoulli(1/2)$ and $B_k \overset{i.i.d.}{\sim} Bernoulli(1/2)$. If we truncate the expansions at depth $d$, then $U_d=\sum_{k=1}^{d}{A_k/2^k}$ and $V_d=\sum_{k=1}^{d}{B_k/2^k}$ are two discrete variables that can take $2^{d}$ possible values. We define the binary variables $\dot{A}_k=2A_k-1$ and $\dot{B}_k=2B_k-1$ to express the interaction between them as their products. We call any products of $A_k$'s and $B_k$'s with at least one $A_k$ and one $B_k$ cross interactions. In other words, cross interactions are defined as the variables of the form $\dot{A}_{k_1}\ldots \dot{A}_{k_r}\dot{B}_{k'_1}\ldots \dot{B}_{k'_t}$ for some $r,t>0$. To explain, we use the following binary integer indexing.  Let $\mathbf{a}$ be a $d$-dimensional binary vector with 1's at $k_1,\ldots,k_r$ and 0's otherwise, and let $\mathbf{b}$ be a $d$-dimensional binary vector with 1's at $k'_1,\ldots,k'_t$ and 0's otherwise. With this notation, the cross interaction $\dot{A}_{k_1}\ldots \dot{A}_{k_r}\dot{B}_{k'_1}\ldots \dot{B}_{k'_t}$ can be written as $\dot{A}_{\mathbf{a}}\dot{B}_{\mathbf{b}}.$

Next, we denote the sum of the observed binary interaction variables by $S_{(\mathbf{ab})}=\sum_{i=1}^n \dot{A}_{\mathbf{a},i} \dot{B}_{\mathbf{b},i}$ with $S_{(\mathbf{00})}=n.$ These statistics are referred to as the symmetry statistics. If $U_d$ and $V_d$ are independent, $(S_{(\mathbf{ab})}+n)/2 \sim Binomial(n,1/2) $ for $\mathbf{a} \neq \mathbf{0}$ and $\mathbf{b}\neq \mathbf{0}$. If marginal distributions are unknown, we can use the empirical CDF transformation and then $(\hat{S}_{(\mathbf{ab})}+n)/4 \sim Hypergeometric(n,n/2,n/2)$ where $\hat{S}_{(\mathbf{ab})}$ is a symmetry statistic with empirical CDF transformation.

If we truncate the expansions at depth $d = d_{max}$, the BET procedure at depth $d_{max}$ can be defined as follows. First, we compute all symmetry statistics with $\mathbf{a} \neq \mathbf{0}$ and $\mathbf{b}\neq \mathbf{0}$ for $d = d_{max}$. For each depth $d = 1, \dots, d_{max}$, we look for the symmetry statistic with the strongest asymmetry and find its $p$-value. Finally, we use Bonferroni adjustment to obtain a $p$-value that considers the family-wise error rate.

The BET has several advantages. The test is minimax optimal under certain regulatory conditions. Moreover, it provides both inferences and clear interpretations. For the BET, rejection of independence implies that there is at least one significant cross interaction. Thus, we can find a potential dependence structure in the sample by investigating the detected cross interaction.

\subsection{Univariate Independence Testing Procedure}

Although the BET shows good performance in many interesting dependency structures, there is room for improvement. In particular, a test based on the sum of squared symmetry statistics provides better power when the sparsity assumption in \citet{zhang2019bet} is violated.

Consider a binary expansion test with specified $d_{max}$. For each depth $d = 1, \dots, d_{max}$, we can find a set of symmetry statistics $S_{(\mathbf{ab})}$. Let $C_d$ be a set of corresponding $\mathbf{ab}$ indices of depth $d$. Since an interaction has different $\mathbf{ab}$ indices for two different $d$, to avoid confusion, we use the $\mathbf{ab}$ of depth $d_{max}$. For example, when $d_{max} = 2$, $C_1 = \{1010\}$ and $C_2 = \{0101, 0110, 0111, 1001, 1010, 1011, 1101, 1110, 1111\}$. The sets $C_d$ have a nested structure. Now, for each depth $d$, we introduce two measures of dependence.
Suppose $X\in\mathbb{R}$ and $Y\in\mathbb{R}$ be two continuous random variables. The population measure of dependence is defined as 
	\begin{eqnarray}
		\mathcal{B}_d(X, Y) = \frac{1}{(2^d-1)^2}\sum_{\mathbf{ab}\in C_d}{E(\dot{A}_{\mathbf{a}}\dot{B}_{\mathbf{b}})^2},
	\end{eqnarray}
	for each depth $d = 1, \dots, d_{max}$.

Let $\{(X_i, Y_i)\}_{i=1}^n$ be a random sample from the joint distribution of $(X, Y)$. The empirical measure of dependence is defined as
	\begin{eqnarray}
		\mathcal{B}_{n,d}[\{(X_i, Y_i)\}_{i=1}^n] = \frac{1}{(2^d-1)^2}\sum_{\mathbf{ab}\in C_d}{\bigg(\frac{S_{(\mathbf{ab})}}{n}\bigg)^2},
	\end{eqnarray}
	for each depth $d = 1, \dots, d_{max}$.
The following theorem lists some properties of $\mathcal{B}_d(X, Y)$ and $\mathcal{B}_{n,d}[\{(X_i, Y_i)\}_{i=1}^n]$.
\begin{theorem}
Suppose the joint distribution of $X$ and $Y$ has a continuous density. The following properties hold:
\begin{enumerate}[\rm(i)]
	\item $\mathcal{B}_d(X, Y) = 0$ if and only if $U_d$ and $V_d$ are independent.
	\item $0\leq\mathcal{B}_d(X, Y)\leq 1$.
	\item $0\leq\mathcal{B}_{n,d}[\{(X_i, Y_i)\}_{i=1}^n]\leq 1$.
	\item $\mathcal{B}_{n,d}[\{(X_i, Y_i)\}_{i=1}^n]\overset{a.s.}{\longrightarrow}\mathcal{B}_d(X, Y)$ as $n\rightarrow \infty$.
	\item If $X$ and $Y$ are independent, then $(2^d-1)^2n\mathcal{B}_{n,d}[\{(X_i, Y_i)\}_{i=1}^n]\overset{d}{\longrightarrow}\chi^2_{(2^d-1)^2}$ as $n\rightarrow \infty$.
\end{enumerate}
\label{THM:uni_property}
\end{theorem}
We define the scaled sum of squared symmetry statistics for each depth $d = 1, \dots, d_{max}$ as
	\begin{eqnarray}
		\xi_{n,d} = \sum_{\mathbf{ab}\in C_d}{\frac{S^2_{(\mathbf{ab})}}{n}}.
	\end{eqnarray}
By this definition, each $\xi_{n,d}$ can be used to detect the dependencies up to depth $d$. Consider a test that rejects $H_0$: ``$X$ and $Y$ are independent'' if at least one of $\xi_{n,d}$'s is greater than $\xi_{n, d, 1-\alpha_d}$, the $1-\alpha_d$ quantile of $\xi_{n,d}$. Then, by Boole's inequality, the upper bound of the type I error is
	\begin{eqnarray}
		Pr(reject~H_0\mid H_0~is~true)\leq \sum_{d=1}^{d_{max}}{\alpha_d}.
	\end{eqnarray}
There are many possible versions of the test based on different choices of the $\alpha_d$'s. Each $\xi_{n,d}$ has a corresponding set of alternatives that it performs well. Therefore, if we have prior information about the dependency, we can choose $\alpha_d$'s in a way that provides optimal power for certain alternatives. When there are no specific prior alternatives, we need a strategy for choosing the $\alpha_d$'s. We remark here that the alternatives in $C_d$ for smaller $d$ are more interpretable than those for larger $d$. From this point of view, we propose an exponentially decaying approach for choice of $\alpha_d$. If we choose $\alpha_d = \alpha \gamma^d/\sum_{d=1}^{d_{max}}{\gamma^d}$  where $0<\gamma\leq 1$ then the upper bound of the significance level is
	\begin{eqnarray}
		Pr(reject~H_0\mid H_0~is~true)\leq \sum_{d=1}^{d_{max}}{\frac{\alpha\gamma^d}{\sum_{d=1}^{d_{max}}{\gamma^d}}}=\alpha,
	\end{eqnarray}
guaranteeing a level $\alpha$ test. A natural choice of $\gamma$ is 1 and the alternatives in each subset, as a group, are equally likely to be detected for $\gamma = 1$;
	\begin{eqnarray}
		Pr(reject~H_0\mid H_0~is~true)\leq \sum_{d=1}^{d_{max}}{\frac{\alpha}{d_{max}}}= \alpha.
	\end{eqnarray}

The power of the proposed test can be improved by a compromise between a distance correlation test and multiple testing over interactions. The binary expansion testing framework loses power from the adverse effect of multiplicity control over depth. This loss of power is particularly severe for linear dependency. See a detailed discussion in Section 1.2 in the supplementary material of \citet{zhang2019bet}. By considering distance correlation combined with the proposed test, we can mitigate this power loss. There is only one interaction in $\xi_{n, 1}$ and it relates to the upper halves of $u$ and $v$ and the lower halves, thus, $\xi_{n, 1}$ represents a linear relationship. Because the above test is composed of multiple hypothesis tests, the test with $\xi_{n, 1}$ can be replaced with the distance correlation test. Because we are using a Bonferroni correction for the critical values, this replacement still maintains the targeted level of the test. We call this approach as ensemble method because it combines two testing methods. The independence test with Pearson's correlation can be also combined with the proposed test. However, we choose the distance correlation test as it improves power in a wider range of cases and it is equivalent to Pearson's correlation under normality. The proposed procedure consists of the following steps:

\begin{enumerate}[\hspace{24pt}(P1)]
\item[\textit{Step 1}]: Fix $\alpha_1, ..., \alpha_{d_{max}}$ with $\sum_{d=1}^{d_{max}}{\alpha_d} = \alpha$.
\item[\textit{Step 2}]: 
Find the $p$-value for the distance correlation test.
\item[\textit{Step 3}]: For each $d = 2, \dots, d_{max}$, compute $\xi_{n, d}$ and its $p$-value.
\item[\textit{Step 4}]: Reject $H_0$ if at least one of the $p$-values is less than respective $\alpha_d$.
\end{enumerate}

To find $p$-value for each depth $d \geq 2$, we can use either a permutation approach or the asymptotic distribution given in theorem \ref{THM:uni_property}, part (v). Now we investigate the behavior of our test in large sample.

\begin{theorem}
Denote the joint distribution of $(U_d, V_d)$ by $\mathbf{P}_{(U_d, V_d)}$ and the bivariate uniform distribution over $\{\frac{0}{2^d}, \dots, \frac{2^d-1}{2^d}\}^2$ by $\mathbf{P}_{0,d}$. For any fixed $0<\delta\leq 1/2$, denote by $\mathcal{H}_{1,d}$ the collection of distributions $\mathbf{P}_{(U_d, V_d)}$ such that $TV(\mathbf{P}_{(U_d, V_d)}, \mathbf{P}_{0,d})\geq \delta$. Consider the testing problem,
\begin{align*}
	H_{0}:\mathbf{P}_{(U_d, V_d)}=\mathbf{P}_{0,d}~~v.s.~~H_{1}:\mathbf{P}_{(U_d, V_d)}\in\mathcal{H}_{1,d}.
\end{align*}
Under $H_1$, each $\xi_{n,d}\rightarrow \infty~as~n\rightarrow \infty$.
\label{THM:uni_consistency}
\end{theorem}

Theorem \ref{THM:uni_consistency} shows that our test statistics, $\xi_{n,d}$'s, go to infinity as sample size increases. Moreover, the distance correlation test is known to be consistent. Therefore, ensemble method is also statistically consistent against the collection of alternatives described in theorem \ref{THM:uni_consistency}.

\subsection{Multivariate Independence Testing Procedure}

Thus far, we have discussed the binary expansion test for univariate random variables. In this section, we develop a generalized independence test for random vectors. The generalization can be made by converting the independence of random vectors into the independence of univariate random variables. The lemma allowing this conversion is stated below.

\begin{lemma}
Let $\mathbf{X}\in \mathbb{R}^p$ and $\mathbf{Y}\in\mathbb{R}^q$ be two random vectors. Then $\mathbf{X}$ and $\mathbf{Y}$ are independent if and only if $\mathbf{s}^T \mathbf{X}$ and $\mathbf{t}^T\mathbf{Y}$ are independent for all $\mathbf{s}\in \mathbb{R}^p, \mathbf{t}\in \mathbb{R}^q$ with $\lVert \mathbf{s} \rVert=1$ and $\lVert \mathbf{t} \rVert=1$.
\label{LEM:Lemma1}
\end{lemma}

This result shows that, to prove independence of random vectors, it is sufficient to consider independence of arbitrary linear combinations of the components. Therefore, the multivariate independence can be tested by checking all possible combinations of $\mathbf{s}$ and $\mathbf{t}$. Because testing all possible combinations cannot be implemented, we consider an approximation of the test by including a finite but reasonably broad number of combinations. Denote hyper unit spheres in $\mathbb{R}^p$ and $\mathbb{R}^q$ by $S_p$ and $S_q$ respectively. Now, for each depth $d$, we propose two measures of dependence for the multivariate setting.

Suppose $\mathbf{X}\in\mathbb{R}^p$ and $\mathbf{Y}\in\mathbb{R}^q$ are two random vectors. For $\mathbf{s}\in S_p, \mathbf{t}\in S_q$, we  define a measure of dependence for the multivariate setting by
	\begin{eqnarray}
		\mathcal{B}_d(\mathbf{X}, \mathbf{Y}) = \frac{1}{c_pc_q}\int_{S_q}{\int_{S_p}{\mathcal{B}_d(\mathbf{s}^T\mathbf{X}, \mathbf{t}^T\mathbf{Y})d\mathbf{s}d\mathbf{t}}},
	\end{eqnarray}
where $c_p=\frac{2\pi^{p/2}}{\Gamma(p/2)}$ and $c_q=\frac{2\pi^{q/2}}{\Gamma(q/2)}$.

Let $\{(\mathbf{X}_i, \mathbf{Y}_i)\}_{i=1}^n$ be a random sample from the joint distribution of $(\mathbf{X}, \mathbf{Y})$. The empirical measure of dependence is defined as
	\begin{eqnarray}
		\mathcal{B}_{n, d}[\{(\mathbf{X}_i, \mathbf{Y}_i)\}_{i=1}^n] = \frac{1}{c_pc_q}\int_{S_q}{\int_{S_p}{\mathcal{B}_{n, d}[\{(\mathbf{s}^T\mathbf{X}_i, \mathbf{t}^T\mathbf{Y}_i)\}_{i=1}^n]d\mathbf{s}d\mathbf{t}}}.
	\end{eqnarray}

The following theorem lists some properties of $\mathcal{B}_d(\mathbf{X}, \mathbf{Y})$ and $\mathcal{B}_{n, d}[\{(\mathbf{X}_i, \mathbf{Y}_i)\}_{i=1}^n]$.

\begin{theorem}
Suppose the joint distribution of $\mathbf{X}$ and $\mathbf{Y}$ has a continuous density. Let $U_d^{\mathbf{s}}$ and $V_d^{\mathbf{t}}$ be truncated binary expansions at depth $d$ of $U^{\mathbf{s}}$ and $V^{\mathbf{t}}$ respectively where $U^{\mathbf{s}}=F_{\mathbf{s}^T\mathbf{X}}(\mathbf{s}^T\mathbf{X})$ and $V^{\mathbf{t}}=F_{\mathbf{t}^T\mathbf{Y}}(\mathbf{t}^T\mathbf{Y})$ for $\mathbf{s}\in S_p, \mathbf{t}\in S_q$. The following properties hold:
\begin{enumerate}[\rm(i)]
	\item $\mathcal{B}_d(\mathbf{X}, \mathbf{Y}) = 0$ if and only if $U_d^{\mathbf{s}}$ and $V_d^{\mathbf{t}}$ are independent for all $\mathbf{s}\in S_p, \mathbf{t}\in S_q$.
	\item $0\leq\mathcal{B}_d(\mathbf{X}, \mathbf{Y})\leq 1$.
	\item $0\leq\mathcal{B}_{n, d}[\{(\mathbf{X}_i, \mathbf{Y}_i)\}_{i=1}^n]\leq 1$.
	\item $\mathcal{B}_{n,d}[\{(\mathbf{X}_i, \mathbf{Y}_i)\}_{i=1}^n]\overset{a.s.}{\longrightarrow}\mathcal{B}_d(\mathbf{X}, \mathbf{Y})~as~n\rightarrow \infty$.
\end{enumerate}
\label{THM:multi_property}
\end{theorem}

Note that $\mathcal{B}_{n,d}[\{(\mathbf{X}_i, \mathbf{Y}_i)\}_{i=1}^n] = E_{\mathbf{S},\mathbf{T}}(\mathcal{B}_{n, d}[\{(\mathbf{S}^T\mathbf{X}_i, \mathbf{T}^T\mathbf{Y}_i)\}_{i=1}^n]\mid\{(\mathbf{X}_i, \mathbf{Y}_i)\}_{i=1}^n)$ where $\mathbf{S}$ and $\mathbf{T}$ follow uniform distributions on $S_p$ and $S_q$ respectively. This expectation can be estimated by
	\begin{eqnarray}
	\widehat{\mathcal{B}}^m_{n, d}[\{(\mathbf{X}_i, \mathbf{Y}_i)\}_{i=1}^n] = \frac{1}{m}\sum_{j=1}^{m}{\mathcal{B}_{n, d}[\{(\mathbf{S}_j^T\mathbf{X}_i, \mathbf{T}_j^T\mathbf{Y}_i)\}_{i=1}^n]},
	\end{eqnarray}
where $\{(\mathbf{S}_j, \mathbf{T}_j)\}_{j=1}^m$ is a random sample generated from uniform distributions on $S_p$ and $S_q$. We call this statistic BERET measure of dependence.

The following theorem shows that the BERET measure of dependence is a consistent estimator of the population measure of dependence.

\begin{theorem}
Suppose the joint distribution of $\mathbf{X}$ and $\mathbf{Y}$ has a bounded continuous density. Then, $\widehat{\mathcal{B}}^m_{n,d}[\{(\mathbf{X}_i, \mathbf{Y}_i)\}_{i=1}^n]\overset{a.s.}{\longrightarrow}\mathcal{B}_d(\mathbf{X}, \mathbf{Y})~as~m,n\rightarrow \infty$.
\label{THM:multi_convergence}
\end{theorem}

Now, to develop an independence test, we define the statistic
	\begin{eqnarray}
		\zeta^m_{n, d} = n(2^d-1)^2	\widehat{\mathcal{B}}^m_{n,d}[\{(\mathbf{X}_i, \mathbf{Y}_i)\}_{i=1}^n],
	\end{eqnarray}
for each depth $d=1, ..., d_{max}$.
By computing $1-\alpha_d$ quantiles of $\zeta^m_{n, d}$, for $d = 1, \dots, d_{max}$, we can consider the test that rejects $H_0:$ ``$\mathbf{X}$ and $\mathbf{Y}$ are independent" if at least one $\zeta^m_{n, d}$, for $d = 1, \dots, d_{max}$, is greater than $\zeta^m_{n, d, 1-\alpha_d}$. If $\sum_{d=1}^{d_{max}}{\alpha_d}\leq \alpha$, this procedure provides a level $\alpha$ test. To put the proposed test into practice, we estimate the asymptotic null distribution by a random permutation method.

For better performance, under possible linear dependency, we combine this procedure with the distance correlation test as above. If the scales of the elements in the random vectors differ greatly, standardization may be helpful to reduce the number of $\mathbf{s}$ and $\mathbf{t}$ to be sampled. We use the normal quantile transformation for standardization. The following procedure summarizes the proposed approach. 

\begin{enumerate}[\hspace{24pt}(P1)]
	\item [\textit{Step 1}]: set $\alpha_1, ..., \alpha_{d_{max}}$ with $\sum_{d=1}^{d_{max}}{\alpha_d} = \alpha$.
	\item [\textit{Step 2}]: standardize marginally each element of the random vectors.
	\item [\textit{Step 3}]: find the $p$-value for the distance correlation test.
	\item [\textit{Step 4}]: fix $m\in \mathbb{N}$ and generate random sample $s_1, \dots, s_m$ and $t_1, \dots, t_m$ from uniform distributions on hyper spheres,  respectivley.
	\item [\textit{Step 5}]: for each $d = 2, \dots, d_{max}$, compute $\zeta^m_{n,d}$ and its $p$-value by the permutation method.
	\item [\textit{Step 6}]: reject $H_0$ if at least one of the $p$-values is less than respective $\alpha_d$.
\end{enumerate}

We refer to this procedure as binary expansion randomized ensemble testing (BERET) due to its two aspects of random projection and ensemble structure. Again we investigate the behavior of our test in large sample. Theorem \ref{THM:multi_consistency} shows that BERET is uniformly consistent against the collection of alternatives in the theorem.

\begin{theorem}
For any fixed $0<\delta\leq 1/2$, denote by $\mathcal{H}^{\mathbf{s},\mathbf{t}}_{1,d}$ the collection of distributions $\mathbf{P}_{(U_d^{\mathbf{s}}, V_d^{\mathbf{t}})}$ such that $TV(\mathbf{P}_{(U_d^{\mathbf{s}}, V_d^{\mathbf{t}})},\mathbf{P}_{0,d})\geq \delta$. Consider the testing problem,
\begin{align*}
	H_{0}:\mathbf{P}_{(U_d^{\mathbf{s}}, V_d^{\mathbf{t}})}=&\mathbf{P}_{0,d}~for~all~\mathbf{s}\in S_p, \mathbf{t}\in S_q\\&v.s.~~H_{1}:\mathbf{P}_{(U_d^{\mathbf{s}}, V_d^{\mathbf{t}})}\in\mathcal{H}^{\mathbf{s},\mathbf{t}}_{1,d}~for~some~\mathbf{s}\in S_p~and~\mathbf{t}\in S_q.
\end{align*}
The following properties hold.
\begin{enumerate}[\rm(i)]
	\item Under $H_1$, $\zeta^m_{n,d}\rightarrow\infty$ as $m,n\rightarrow\infty$.
	\item Rejection probability of the permutation test is bounded by $\alpha$ under $H_0$ and converges to 1 under $H_1$ as $m, n \rightarrow \infty$ if $d_{max} \geq d$.
\end{enumerate}
\label{THM:multi_consistency}
\end{theorem}

The BERET has the following advantages. First, the method achieves robust power by a compromise between the distance correlation test and multiple testing over interactions. The power loss due to multiplicity control over the depth also exists in the multivariate case. By considering the distance correlation result together with the proposed measure of dependence with $d \ge 2$, we can improve power over a wide range of plausible dependencies.

\begin{figure}[h]
\centering
\includegraphics[scale = 0.5]{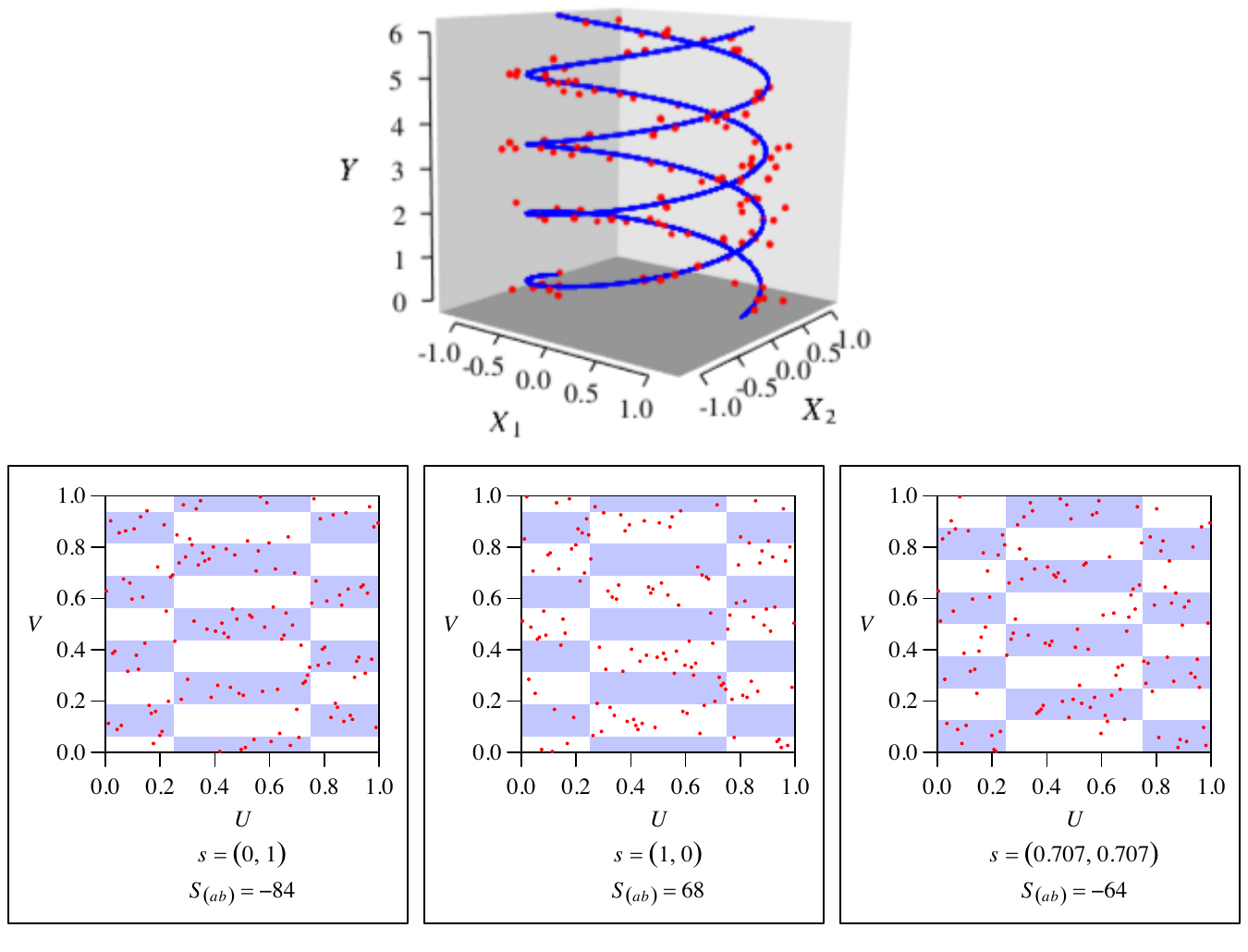}
\caption{Analysis for double helix dependency. The plot in the first row shows a sample with double helix dependency between a random vector $(X_1~X_2)^T$ and a random variable $Y$ with $n = 128$. The blue line is added to illustrate the true relationship. The three plots in the second row show the linear combinations of $X_1$ and $X_2$ with the strong asymmetries and the corresponding symmetry statistics $(S_{(\mathbf{ab})})$. Positive regions ($\dot{A}_{\mathbf{a}}\dot{B}_{\mathbf{b}} = 1$) are in white and negative regions ($\dot{A}_{\mathbf{a}}\dot{B}_{\mathbf{b}} = -1$) are in blue.}
\label{FIG:interpretability}
\end{figure}

The second benefit of our method is clear interpretability. The issue of interpretability is particularly important in evaluating multivariate relationships. However, most multivariate independence tests provide only the results of the tests with no information on potential dependence structures in the sample. Although the canonical correlation test provides some related information, it shows poor power in nonlinear relationships relative to the proposed method. In contrast, when the proposed test rejects independence, the $\mathbf{s}$ and $\mathbf{t}$ vectors indicate the linear combinations of the vectors that have strong dependencies. Using these vectors, we can detect the possible dependence structures in the sample. See Fig. \ref{FIG:interpretability} for a three dimensional double helix structure example for illustration, in which white positive regions and blue negative regions of interaction provide the interpretation of global dependency. It can be seen that the double helix structure is detected by three linear combinations. More interesting interpretation examples are provided in Section \ref{data}.

Lastly, our method provides useful exploratory information for model selection. A small entry in the unit vector $\mathbf{s}$ or $\mathbf{t}$ may indicate that the corresponding variable may not be related to the other random vector. See data examples in Section \ref{data} for details.

\section{Simulation Studies}\label{simulation}

\subsection{Univariate Independence}
For comparison, we consider the Hoeffding's D test, the distance correlation test, the mutual information test (MINTav), Fisher's exact scanning method, and the maximum binary expansion test. We use sample size $n = 128$ as a moderate sample size for power comparison. We set the level of the tests to be 0.1 and simulate each scenario 1,000 times. We adopt $d_{max} = 4$ because this depth provides a good approximation to the true distribution. See a detailed discussion in Section 4.5 in \citet{zhang2019bet}. The p-values of the proposed method are calculated using the asymptotic distribution of theorem \ref{THM:uni_property}, part (v). We verified that the p-value under the null hypothesis was controlled at level 0.1.

\begin{figure}[!tbph]
\centering
\includegraphics[scale = 0.7]{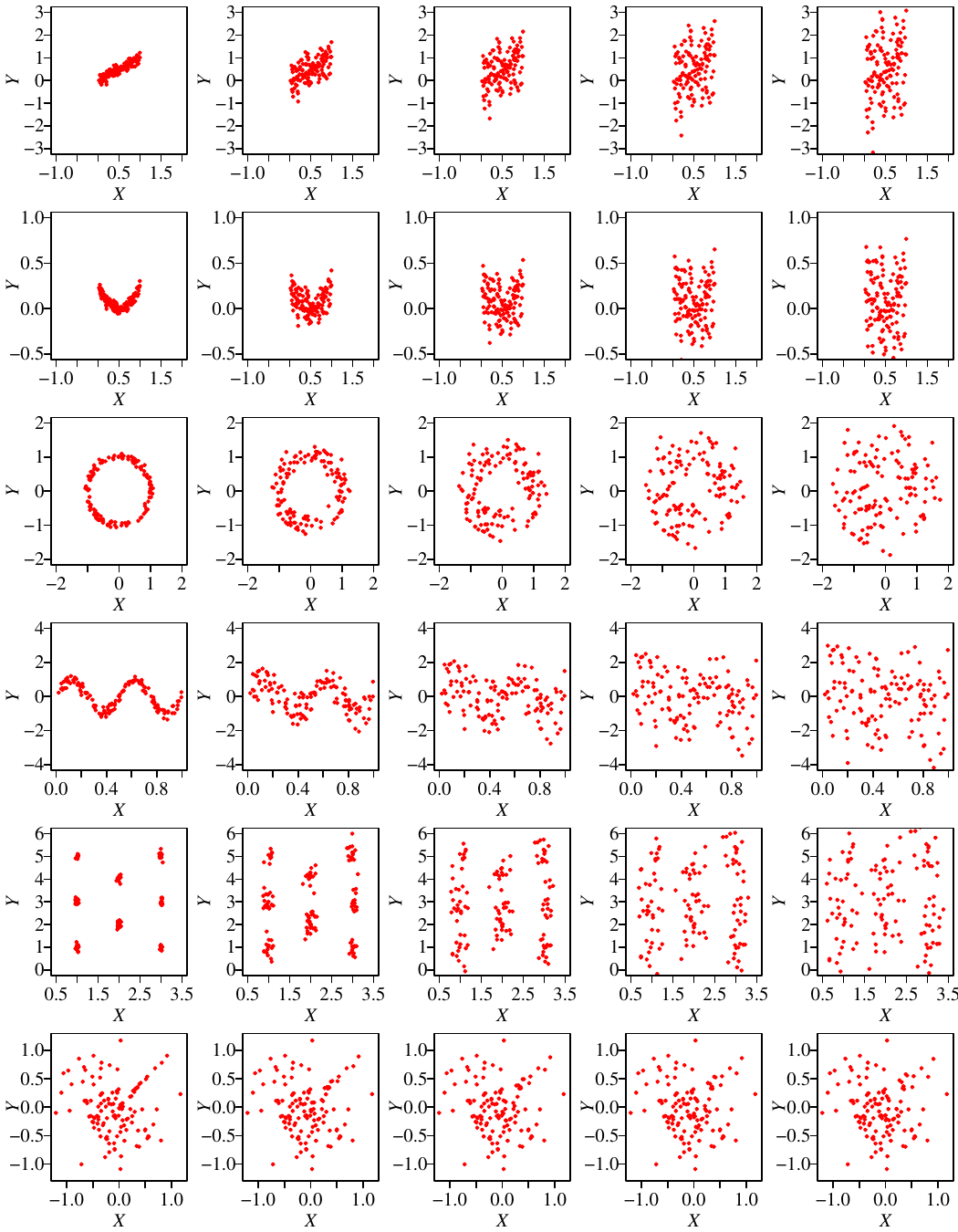}
\caption[Visualization of scenarios in univariate independence test]{The observations are generated from the six scenarios in Table 1 with $n = 128$ and noise levels $l = 1, 3, 5, 7,$ and $9$. Each row shows the scatter plots of linear, parabolic, circular, sine, checkerboard, and local dependency scenarios in order.}
\label{FIG:scenario_example}
\end{figure}

We compare the power of the above methods over familiar dependence structures such as linear, parabolic, circular, sine, checkerboard and local relationship described in \citet{zhang2019bet}. At each noise level $l=1,\ldots,10$, $\epsilon, \epsilon', \epsilon''$ are independent $\mathcal{N}(0,(l/40)^2)$ random variables. $U$ follows the standard uniform distribution. $\vartheta$ is a $U[-\pi,\pi]$ random variable. $W$, $V_1$, and $V_2$ follow $multi$-$Bern(\{1,2,3\}, (1/3,1/3,1/3))$, $Bern(\{2,4\},(1/2,1/2))$, and $multi$-$Bern(\{1,3,5\}, (1/3,1/3,1/3))$ respectively. $G_1$, $G_2$ are generated from $\mathcal{N}(0,1/4)$. Table \ref{table:univariate_scenarios} summarizes the details of the setting. Some graphical descriptions of the scenarios are given in Fig. \ref{FIG:scenario_example}.

\begin{table*}[h]
\caption{Simulation scenarios for univariate independence test}
\label{table:univariate_scenarios}
\centering
\scriptsize
\begin{tabular}{@{}lll@{}}
\hline
Scenario & Generation of $X$ & Generation of $Y$\\
\hline
Linear & $X=U~~~~~$ & $ Y=X+6\epsilon$ \\
Parabolic & $X=U~~~~~$ & $  Y=(X-0.5)^2+1.5\epsilon$\\
Circular & $X=\cos \vartheta+2\epsilon~~~~~$ & $  Y=\sin\vartheta+2\epsilon'$\\
Sine & $X=U~~~~~$ & $  Y=\sin(4 \pi X)+8 \epsilon$\\
Checkerboard & $X=W+\epsilon~~~~~$ & $  Y=\left\{\begin{array}{ll}
V_1 +4\epsilon' & \text{ if }  W=2\\
V_2 +4\epsilon'' & \text{ otherwise}
\end{array}\right.$\\
Local & $X=G_1~~~~~$ & $  Y=\left\{\begin{array}{ll}
X+\epsilon &  \text{ if } 0 \le G_1 \le 1 \text{ and } 0 \le G_2 \le 1\\
G_2 & \text{ otherwise}
\end{array}\right.$\\
\hline
\end{tabular}
\end{table*}

\begin{figure}[!tbph]
\centering
\includegraphics[scale = 0.68]{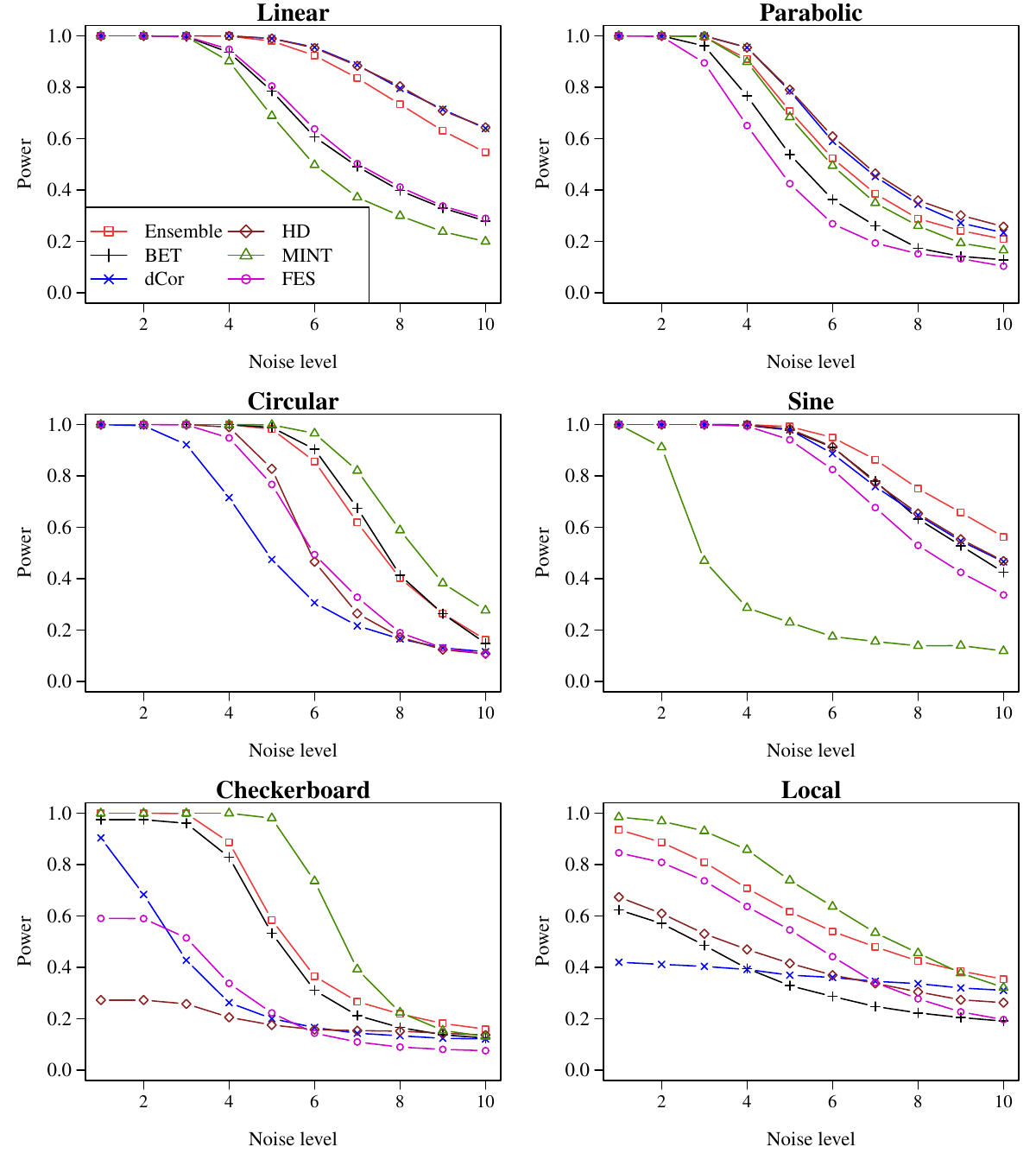}
\caption[Comparison of powers in univariate dependencies]{Comparison of powers from six tests of independence: the binary expansion randomized ensemble test with $d_{max}=4$ (red square), the maximum binary expansion test with $d_{max}=4$ (black plus sign), the distance correlation test (blue cross), Hoeffding's D(brown diamond), the mutual information test  (green triangle), and Fisher exact scanning (purple circle).}
\label{FIG:univariate_results}
\end{figure} 

Fig. \ref{FIG:univariate_results} shows the performance of the six methods. There are two points to notice. First, except for the proposed test, all the other methods show the lowest power in at least one scenario. The ensemble approach and the BET show similar powers across the scenarios except for the linear and local dependency. The ensemble approach considerably improves power in the linear and local dependency scenarios. As discussed previously, the ensemble approach utilizes the information on dependence remaining in the symmetry statistics that is not reflected in calculation of the maximum binary expansion testing. Therefore, small asymmetries in many symmetry statistics can be combined to provide a significant result in the ensemble approach when the sparsity assumption is violated. This result is related to the second finding that the ensemble approach outperforms Fisher's exact scanning in both global and local dependence structures. \citet{zhang2019bet} reported that maximum binary expansion testing provides better power for global dependence structures, whereas Fisher's exact scanning performs better for local dependence structures. The simulation results suggest that the ensemble approach works better than Fisher's exact scanning even in the local dependency scenario. 

\subsection{Multivariate Independence}

Although the proposed method can be applied to arbitrary $p$ and $q$, we choose $p = 2$ and $q = 1$ for better illustration. We compare the proposed method with the distance correlation test, the Heller-Heller-Gorfine test, the $d$-variable Hilbert-Schmidt independence criterion, and the mutual information test (MINTav). We again use sample size $n = 128$. We set the level of the tests to be 0.1 and simulate each scenario 1,000 times. For our method, we adopt $m = 30$ because there is no considerable difference in performance compared with larger $m$'s such as $m = 360$. We also use a permutation method with 1,000 replicates to calculate the p-values of the proposed approach. We verified that the p-value under the null hypothesis was controlled at level 0.1.

\begin{figure}[!tbph]
\centering
\includegraphics[scale = 0.68]{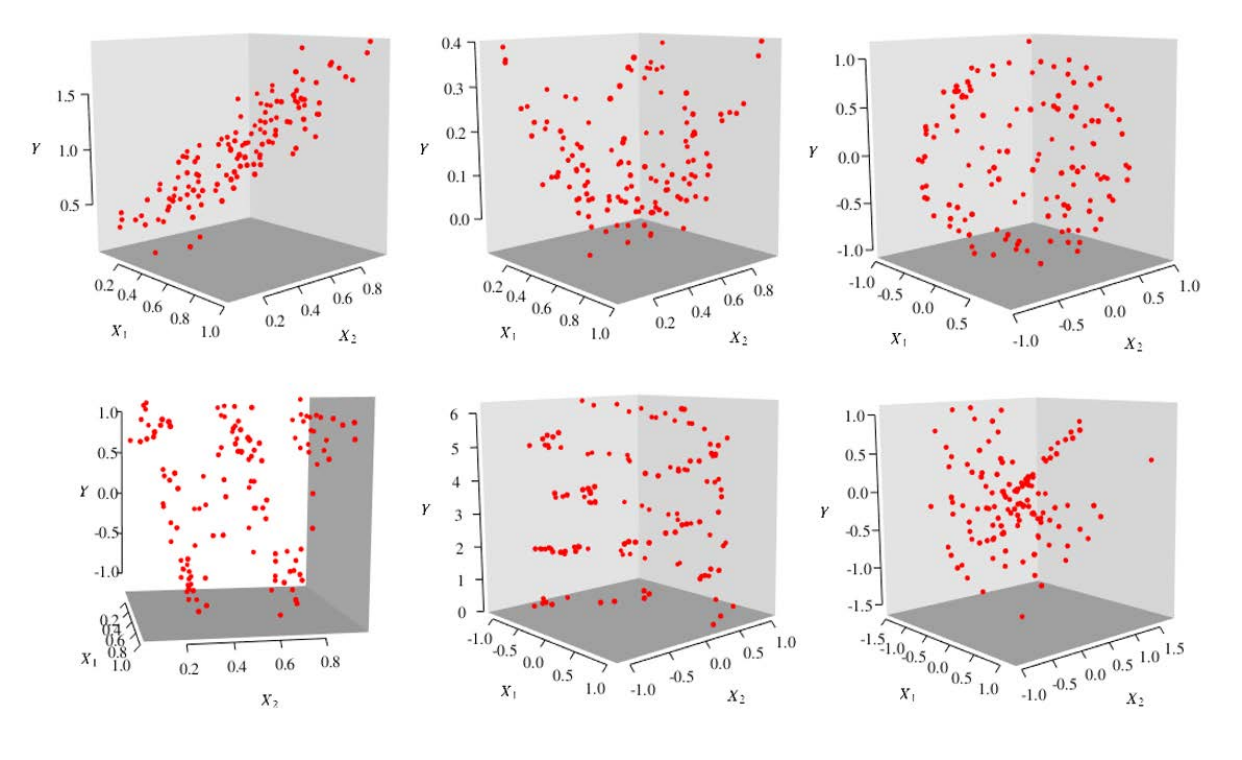}
\caption[Visualization of scenarios in multivariate independence test]{The observations are generated from the six scenarios in Table 2 with $n = 128$ and noise level $l = 1$. Each plot shows the scatterplots of linear (upper left), parabolic (upper center), spherical (upper right), sine (bottom left), double helix (bottom center), and local dependency scenarios (bottom right).}
\label{FIG:scenario_example_3d}
\end{figure}

We compare the power of the methods over dependence structures such as linear, parabolic, spherical, sine, and local dependence structures. These scenarios are generalized from the univariate dependence simulations. In addition, we include an additional interesting relationship, the double helix structure. At each noise level $l=1,\ldots,10$, $\epsilon, \epsilon', \epsilon''$ are independent $\mathcal{N}(0,(l/40)^2)$ random variables. $U_1, U_2$ follow the standard uniform distribution. $\vartheta$ follows $\mathrm{U}[0,4\pi]$. $G_1, G_2, G_3$ are independent $\mathcal{N}(0,1/4)$ random variables. $I$ follow the Rademacher distribution. Table \ref{table:multivariate_scenarios} summarizes the details of the setting. These three dimensional scenarios are visually displayed in Fig. \ref{FIG:scenario_example_3d}. 

\begin{table*}[!tbph]
\caption{Simulation scenarios for multivariate independence testing}
\label{table:multivariate_scenarios}
\centering
\scriptsize
\begin{tabular}{@{}lll@{}}
\hline
Scenario & Generation of $X$ & Generation of $Y$\\
\hline
 Linear & $\mathbf{X} = \begin{pmatrix} U_1 \\ U_2 \end{pmatrix}$ & $\mathbf{Y} = X_1 + X_2 + 7\epsilon$   \\ 
 Parabolic & $\mathbf{X} = \begin{pmatrix} U_1 \\ U_2 \end{pmatrix}$ &$\mathbf{Y} = (X_1-0.5)^2 + (X_2-0.5)^2 + 1.5\epsilon$   \\ 
 Spherical & $\mathbf{X} = \begin{pmatrix} \frac{G_1}{\sqrt{G_1^2+G_2^2+G_3^2}} \\ \frac{G_2}{\sqrt{G_1^2+G_2^2+G_3^2}} \end{pmatrix}$ & $\mathbf{Y} = \frac{G_3}{\sqrt{G_1^2+G_2^2+G_3^2}} + 3 \epsilon$  \\ 
 Sine & $\mathbf{X} = \begin{pmatrix} U_1 \\ U_2 \end{pmatrix}$ & $\mathbf{Y} = \sin{(5\pi X_1)}+ 4\epsilon$   \\ 
 Double helix & $\mathbf{X} = \begin{pmatrix} Icos{(\vartheta) + 1.5\epsilon} \\ Isin{(\vartheta) + 1.5\epsilon'} \end{pmatrix}$ & $\mathbf{Y} = \frac{\vartheta}{2} + 2\epsilon''$  \\ 
 Local & $\mathbf{X} = \begin{pmatrix} G_1 \\ G_2 \end{pmatrix}$ & $\mathbf{Y} = \begin{cases}
    \frac{X_1}{\sqrt{2}} + \frac{X_2}{\sqrt{2}} + \frac{\epsilon}{2}, & \text{if $0\leq G_1 + G_2\leq 2$ and $0\leq G_3\leq 1$}.\\
    G_3, & \text{otherwise}.
  \end{cases}$ \\
\hline
\end{tabular}
\end{table*}

\begin{figure}[!tbph]
\centering
\includegraphics[scale = 0.65]{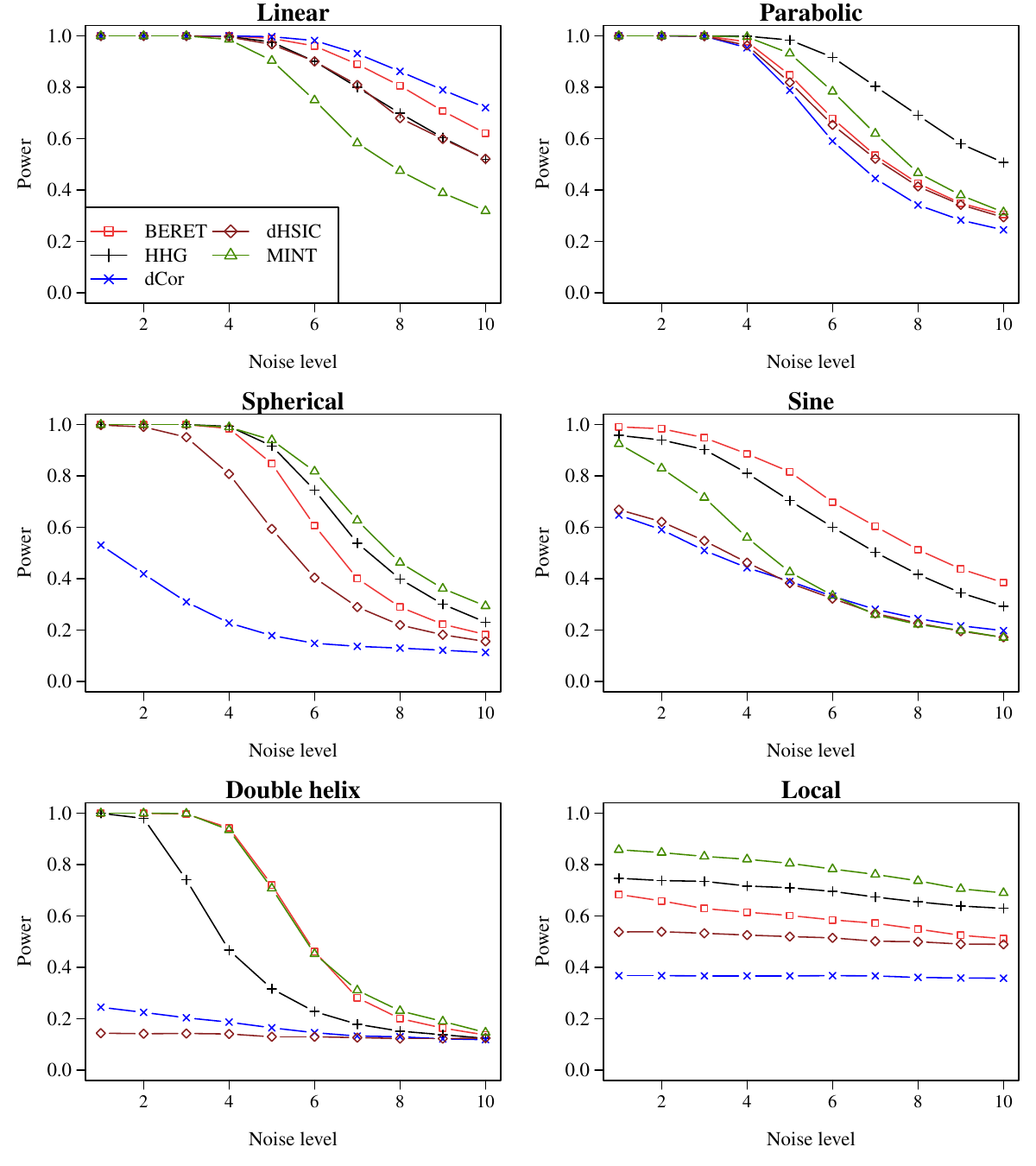}
\caption[Comparison of powers in multivariate dependencies]{Comparison of powers from five tests of independence: the binary expansion randomized ensemble test with $d_{max}=4$ (red square), the Heller-Heller-Gorfine test (black plus sign), the distance correlation test (blue cross), the $d$-variable Hilbert-Schmidt independence criterion (brown diamond), and the mutual information test (green triangle).}
\label{FIG:multivariate_results}
\end{figure}

Fig. \ref{FIG:multivariate_results} shows the simulation results. The BERET provides the best power in more complex dependency structures such as sine and double helix dependency. It outperforms the distance correlation test and the $d$-variable Hilbert-Schmidt independence criterion in at least five scenarios compared with each testing method separately. Moreover, our method provides stable results across the scenarios considered. It ranks at least third place in all scenarios. The mutual information test performs best in the highest number of scenarios. In linear and sine relationships, however, there is significant loss of power in the mutual information test compared with the proposed method. A point to notice is that our method provides additional insight. Other methods only provide test results of independence, but our method provides potential dependence structures as well as test results. The simulation results show that BERET provides competitive performance while providing a much clearer interpretation.

\section{Data Examples} \label{data}

\subsection{Life Expectancy}

\begin{figure}[!tbph]
\centering
\includegraphics[scale = 0.85]{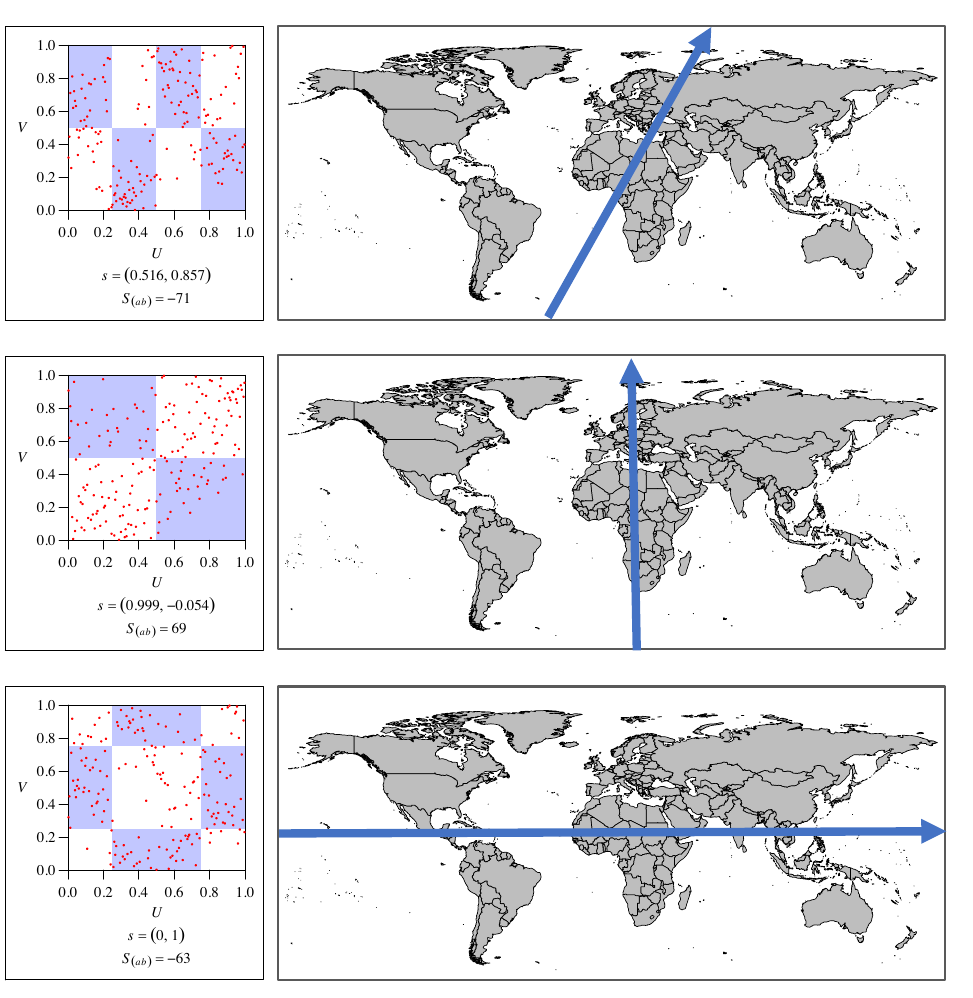}
\caption[Detected dependency structures of life expectancy]{The three strongest dependency structures between geographic location and life expectancy. They also present the corresponding values of the symmetry statistics $(S_{(\mathbf{ab})})$ and the coefficients of the linear combination $(\mathbf{s})$ of $X_1$ and $X_2$. The blue arrows in the world maps represent the horizontal axes in the scatterplots.}
\label{FIG:life_expectancy1}
\end{figure}

We use the proposed method to test independence between geographic location and life expectancy and compare its performance with the performance of other methods, i.e., the distance correlation test (dCor), the Heller-Heller-Gorfine test (HHG), the mutual information test (MINT), and the canonical correlation test (CC). We include the canonical correlation test because it also provides some insight on dependence structure as does the proposed method. For the proposed method, we set $d_{max} = 4$ and $m = 30$. The p-value of the test is calculated by a permutation method with 1,000 replicates. The dataset is obtained from the life expectancy report released by the World Health Organization in 2016. The dataset includes males and females and total life expectancy of 189 countries and special administrative regions estimated in 2015. The latitudes ($X_1$), longitudes ($X_2$), and total life expectancies ($Y$) are used in the analysis. Table \ref{table:life_expectancy1} presents the testing results for the five different methods. All five tests provide p-values close to 0, indicating a significant dependence between geographic location and life expectancy. To identify the dependence structure, we investigate the symmetry statistics. Fig. \ref{FIG:life_expectancy1} shows the three largest symmetry statistics and the corresponding $\mathbf{s}$'s.

\begin{table*}[h]
\caption{p-values from five tests of independence for life expectancy}
\label{table:life_expectancy1}
\centering
\begin{tabular}{lrrrrr}
	\hline
	& \multicolumn{1}{c}{BERET} & \multicolumn{1}{c}{dCor} & \multicolumn{1}{c}{HHG} & \multicolumn{1}{c}{MINT} & \multicolumn{1}{c}{CC} \\
	\hline
Original sample & $<$0.0001 & $<$0.0001 & 0.0010 & 0.0010 & $<$0.0001 \\
With noise & $<$0.0001 & $<$0.0001 & 0.0010 & 0.0010 & $<$0.0001 \\
$n = 64$ & 0.0020 & 0.0014 & 0.0050 & 0.0010 & 0.0158 \\
$n = 32$ & 0.0806 & 0.0652 & 0.0877 & 0.0177 & 0.0995 \\
\hline
\end{tabular}
\end{table*}

The most asymmetric result is shown in the first row. It is $\dot{A}_2\dot{B}_1$ with $\mathbf{s} = (0.516,~0.857)^T$. The horizontal axis is the empirical cumulative distribution function transformation of $0.516 X_1 + 0.857 X_2$, wherein a smaller value implies the country is located in the southwest and a larger value implies the northeast. There are four different groups, from the first one in the upper left to the fourth group in the lower right. Each blue cell represents a specific region, America, Africa, Europe and Asia from left to right. The countries in America and Europe show higher life expectancy than countries in Africa and Asia. The four points in the top right corner are Hong Kong, Japan, Macau and South Korea. They can be interpreted as potential outliers distinct from the global pattern. 

The second row shows that there is a positive relationship between latitude and life expectancy. That is, the countries in North America, Europe and Northeast Asia have higher life expectancy than countries in Africa, South America and the other parts of Asia. The last row shows that a circular dependency can exist, which indicates that countries in the America and Asia have a medium life expectancy, whereas countries around the prime meridian have different life expectancies, higher in Europe and lower in Africa. These findings prove clearly that our method detects the dependence structures between geographic location and life expectancy.
 
The canonical correlation analysis also can be used to find information on dependence structure. The canonical correlation is 0.43, and it is calculated using $0.991X_1 - 0.137X_2$ and $Y$. The coefficients of $X_1$ and $X_2$ are similar to the elements of $\mathbf{s}$ in the result of the proposed method in the second row. However, canonical correlation provides information only on the linear dependence structure, whereas our method provides richer information by considering various nonlinear dependence structures.

\begin{figure}[!tbph]
\centering
\includegraphics[scale = 0.7]{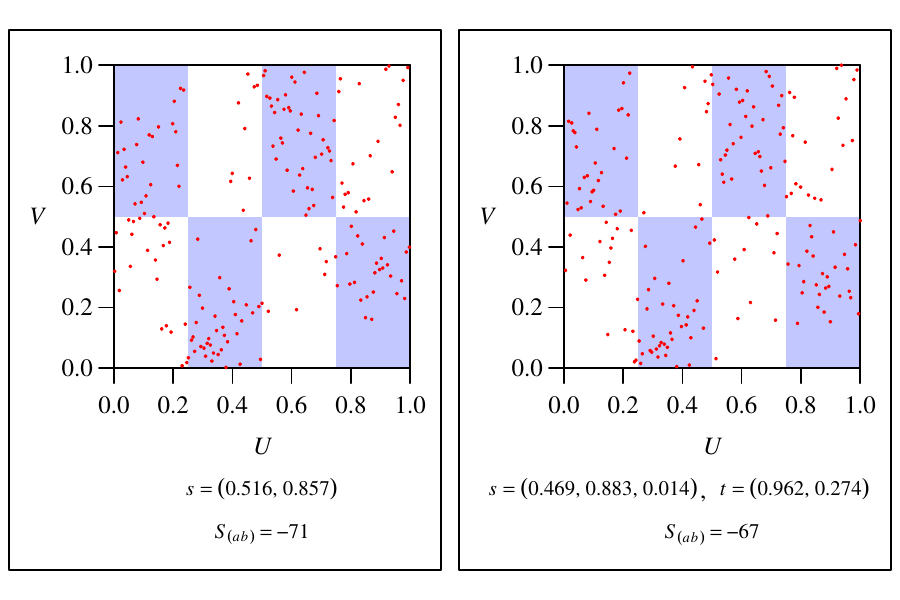}
\caption[Detected dependency structures of life expectancy with noise]{The left plot shows the strongest dependency structures between geographic location and life expectancy. The right plot shows the corresponding result after inclusion of noise variables. The third entry of the vector $\mathbf{s}$ and the second entry of the vector $\mathbf{t}$ in the right plot are the coefficients of noise variables. The plots also present the values of the symmetry statistics $(S_{(\mathbf{ab})})$ and the coefficients in the linear combinations $\mathbf{s}$ and $\mathbf{t}$.}
\label{FIG:life_expectancy_noise}
\end{figure}

Now we add two randomly sampled noise variables to each of $\mathbf{X}$ and $\mathbf{Y}$ to evaluate performance of the proposed method in the presence of noisy data. Each noise variable is chosen from a standard normal distribution. Table \ref{table:life_expectancy1} presents the p-values for the five different methods. The p-values of all five methods are robust. Fig. \ref{FIG:life_expectancy_noise} shows the previously detected strongest dependence structure without noise variables along with an example of the corresponding result with noise variables. The coefficients of the variables are stable and the coefficients of the noise variables are small.

To evaluate small sample performance, we randomly select small subsamples from the original sample. To account for randomness, we calculate the average p-values from 100 random subsamples. Table \ref{table:life_expectancy1} displays the average p-values of all five methods. All p-values are similarly affected, and all of them increase as sample size decreases.

In summary, all five methods detect the dependence with very small p-values. In addition, the BERET detects three interesting dependence structures which can be explained by known global features. The effects of inclusion of noise variables and reduction in sample size are similar for the five methods.

\subsection{Mortality Rate}
The second case is the relationship between mortality rate, birth rate and income level. We use the Central Intelligence Agency's world fact data, estimated in 2018. The dataset includes income level ($X_1$), birth rate ($X_2$),  and mortality rate ($Y$) of 225 countries and special administrative regions. Table \ref{table:mortality_rate1} presents the calculated p-values.

\begin{table*}[h]
\caption{p-values from five tests of independence for mortality rate}
\label{table:mortality_rate1}
\centering
\begin{tabular}{lrrrrr}
	\hline
	& \multicolumn{1}{c}{BERET} & \multicolumn{1}{c}{dCor} & \multicolumn{1}{c}{HHG} & \multicolumn{1}{c}{MINT} & \multicolumn{1}{c}{CC} \\
	\hline
Original sample & 0.0040 & 0.0050 & 0.0010 & 0.3077 & 0.4303 \\
With noise & 0.0231 & 0.0213 & 0.0020 & 0.3671 & 0.5272 \\
$n = 64$ & 0.0287 & 0.3288 & 0.1768 & 0.4998 & 0.4350 \\
$n = 32$ & 0.2812 & 0.4631 & 0.3304 & 0.4778 & 0.4359 \\
\hline
\end{tabular}
\end{table*}

\begin{figure}[h]
\centering
\includegraphics[scale = 0.54]{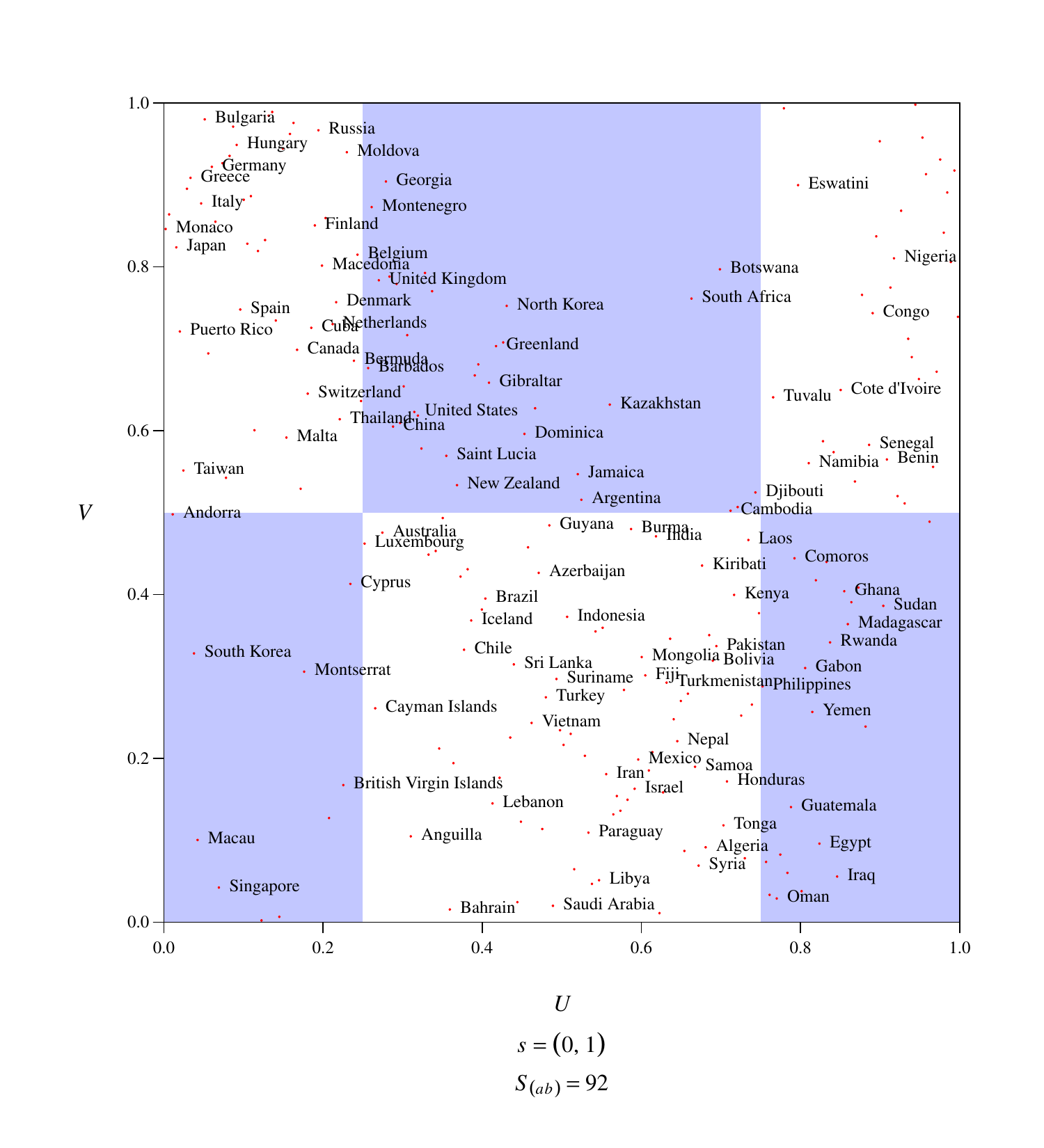}
\caption[Detected dependency structure of mortality rate]{The plot shows the strongest dependency structure between birth rate, income level and mortality rate. It also presents the corresponding value of the symmetry statistic $(S_{(\mathbf{ab})})$ and the coefficients of the linear combination $(\mathbf{s})$ of $X_1$ and $X_2$.}
\label{FIG:mortality_rate1}
\end{figure}

Once again, the proposed method and two other methods provide p-values close to 0, which rejects the null hypothesis, whereas the mutual information test and canonical correlation fail to reject it. The poor performance of canonical correlation can be explained by investigating the results of our method. The strongest asymmetry is given in Fig. \ref{FIG:mortality_rate1}, which shows a strong quadratic relationship. This relationship explains the failure of canonical correlation to work for the data we use here. Although the canonical correlation test provides both inference and information on dependence structure, it performs poorly in nonlinear dependency settings.

For explanation of the observed quadratic relationship one must point to two conflicting phenomena. The first one is that in developed countries the birth rates are low, but the mortality rates are high due to population aging. In developing countries, however, the birth rates are high from lack of family planning and the mortality rates are also high due to insufficient public health. Thus, mortality rates are high in countries with low or high birth rates.

\begin{figure}[h]
\centering
\includegraphics[scale = 0.7]{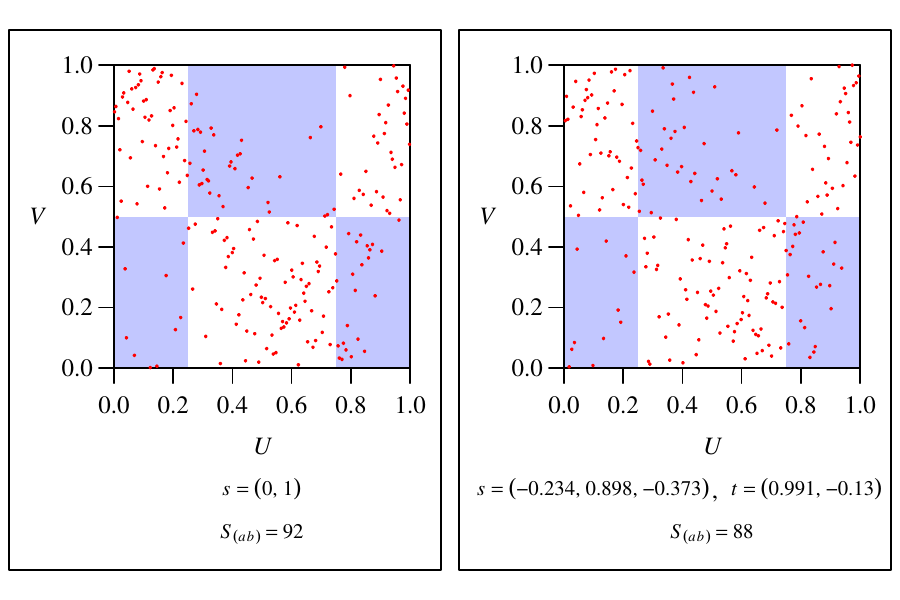}
\caption[Detected dependency structures of mortality rate with noise]{The plots show the strongest dependency structure between birth rate, income level and mortality rate with and without noise variables. The third entry of the vector $\mathbf{s}$ and the second entry of the vector $\mathbf{t}$ in the right plot are the coefficients of noise variables. The plots also present the values of the symmetry statistics $(S_{(\mathbf{ab})})$ and the coefficients in the linear combinations $\mathbf{s}$ and $\mathbf{t}$.}
\label{FIG:mortality_rate_noise}
\end{figure}

We investigate the effect of adding noise variables in the same manner as above. The p-values of the five different methods are represented in Table \ref{table:mortality_rate1}, which shows that all five methods are similarly affected by the noise variables. The strongest dependence structures detected by our method are presented in Fig. \ref{FIG:mortality_rate_noise}. The coefficients of the noise variables are relatively small, and the same dependence structures are identified as before.

For small sample performance comparison, we randomly select small subsamples as we did previously. Table \ref{table:mortality_rate1} lists the average p-values. The result shows that the p-values of the distance correlation test and the Heller-Heller-Gorfine test are significantly more affected on average than the p-values of the proposed method.

Once again, the BERET detects interesting structure that can be explained by widely recognized relationships between mortality rate and birth rate. It provides stable performance even when there are noise variables or when the sample size is small, whereas other methods can be significantly affected by a reduction in sample size.

\subsection{House Price}

The third data example is the market historical dataset of real estate from the University of California, Irvine machine learning repository. The data includes 414 transactions from the Xindan district of Taipei between August 2012 and July 2013. We use these data to detect the relationship between geographic location and house price. The p-values of the five methods are presented in Table \ref{table:house_price1}. 

\begin{table*}[h]
\caption{p-values from five tests of independence for house price}
\label{table:house_price1}
\centering
\begin{tabular}{lrrrrr}
	\hline
	& \multicolumn{1}{c}{BERET} & \multicolumn{1}{c}{dCor} & \multicolumn{1}{c}{HHG} & \multicolumn{1}{c}{MINT} & \multicolumn{1}{c}{CC} \\
	\hline
Original sample & $<$0.0001 & $<$0.0001 & 0.0010 & 0.6204 & $<$0.0001 \\
With noise & 0.0040 & 0.4910 & 0.4961 & 0.5206 & $<$0.0001 \\
$n = 64$ & $<$0.0001 & $<$0.0001 & 0.0010 & 0.2809 & $<$0.0001 \\
$n = 32$ & $<$0.0001 & $<$0.0001 & 0.0016 & 0.2660 & 0.0067 \\
\hline
\end{tabular}
\end{table*}

\begin{figure}[!tbph]
\centering
\includegraphics[scale = 0.93]{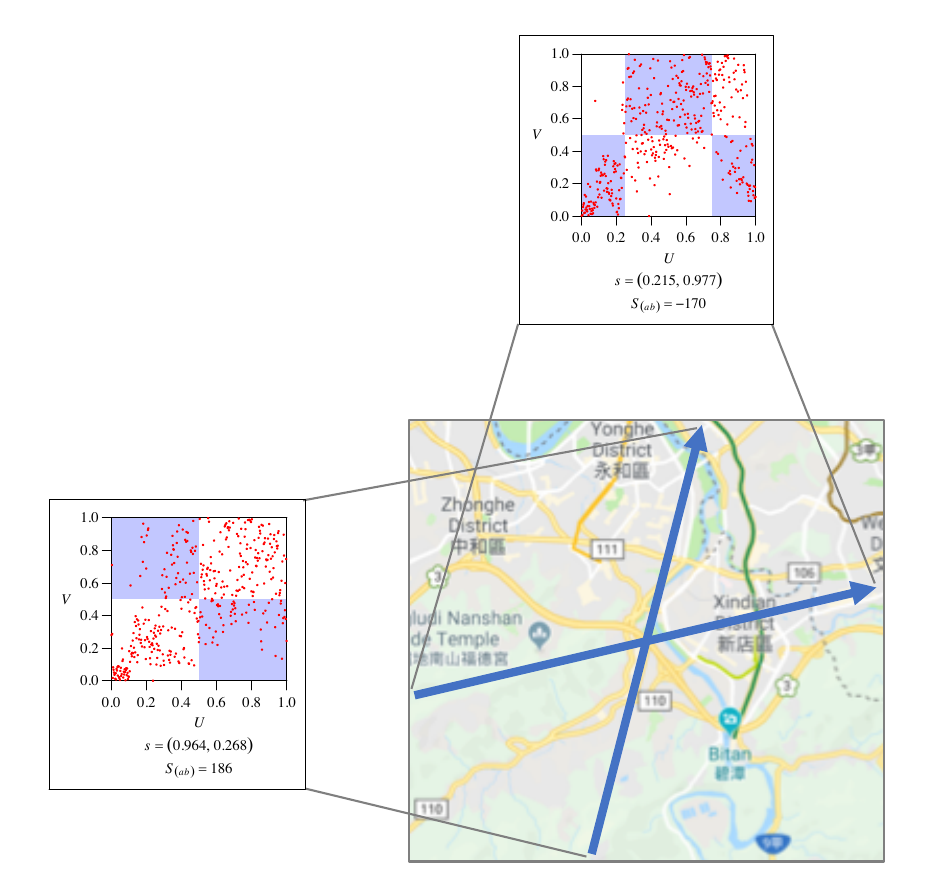}
\caption[Detected dependency structures of house price]{The plots show the two strongest dependency structures between geographic location and house price. The plots also present the values of the symmetry statistics $(S_{(\mathbf{ab})})$ and the coefficients in the linear combinations $\mathbf{s}$ and $\mathbf{t}$. The blue arrows in the map represent the horizontal axes in the scatterplots.}\label{FIG:house_price1}
\end{figure}

All methods except the mutual information test provide p-values close to 0, which is consistent with the commonly assumed relationship between location and house price in a city. The mutual information test fails to reject the independence. Fig. \ref{FIG:house_price1} presents the two strongest dependencies identified by the proposed method.

The symmetry statistic with the strongest asymmetry is $\dot{A}_1\dot{B}_1$, which means that there may be a linear relationship between geographic location and house price. The corresponding $\mathbf{s}$ for the horizontal axis is $(0.964, 0.268)$. That is, houses have higher values in the north and lower values in the south. It is because the central part of Taipei is above the Xindan district. The symmetry statistic with the second strongest asymmetry is $\dot{A}_1\dot{A}_2\dot{B}_1$. The corresponding $\mathbf{s}$ for the horizontal axis is $(0.215,0.977)^T$. That is, house prices are high at the center of the district, where two main roads intersect, and prices fall towards the periphery. These results accord closely with the general characteristics of real estate prices in a city. Therefore, we can conclude that the proposed method properly detects the relationships between house price and geographic location.

\begin{figure}[!tbph]
\centering
\includegraphics[scale = 0.7]{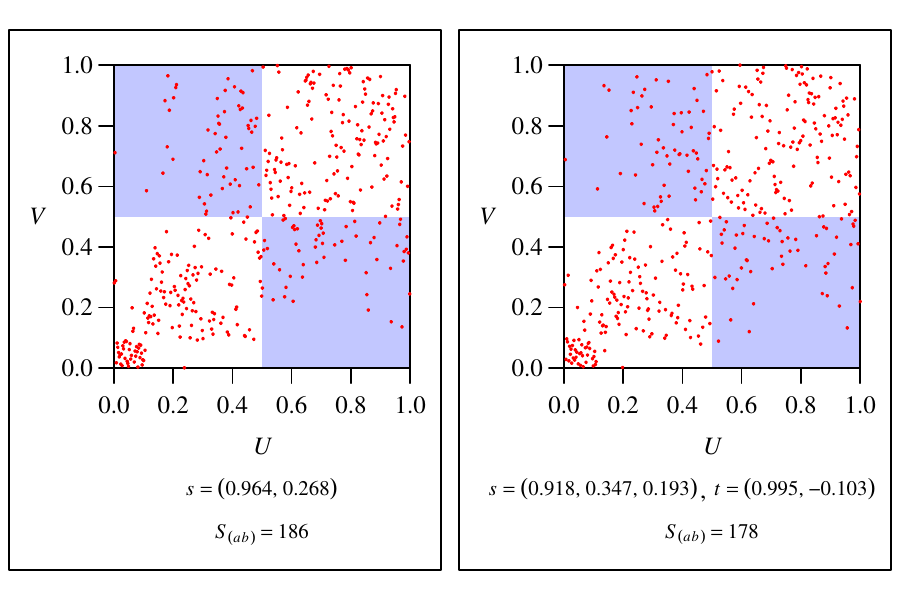}
\caption[Detected dependency structures of house price with noise]{The plots show the strongest dependency structure between geographic location and house price detected by the proposed method without noise variables and the corresponding result with noise variables. The third entry of the vector $\mathbf{s}$ and the second entry of the vector $\mathbf{t}$ in the right plot are the coefficients of the noise variables. The plots also present the values of the symmetry statistics $(S_{(\mathbf{ab})})$ and the coefficients in the linear combinations $\mathbf{s}$ of $X_1$,  $X_2$, and $X_3$ and $\mathbf{t}$ of $Y_1$ and $Y_2$.}
\label{FIG:house_price_noise}
\end{figure}

We add two randomly sampled noise variables to each $\mathbf{X}$ and $\mathbf{Y}$ as before. The resulting p-values of the five different methods are represented in Table \ref{table:house_price1}. The results of all methods except our method are significantly affected by the noise variables. The possible dependence structures detected by our method are presented in Fig. \ref{FIG:house_price_noise}. The figure indicates that the same dependence structures are detected and the coefficients of the noise variables are relatively small as before. The average p-values from the small random subsamples are also given in Table \ref{table:house_price1}. The result shows that there is little difference in power among the four significant methods.

\section{Conclusion} \label{conclusion}

Detection of dependence in a distribution-free setting is an important problem in statistics. Existing methods may have challenges with detecting complicated dependence structures. The distance correlation test, for example, does not detect circular dependency well, whereas it provides good powers in linear, parabolic, and sine settings in simulation studies. The binary expansion testing procedure in \citet{zhang2019bet} suggests a novel way to solve this problem. However, it is limited to the independence test of two random variables and there is room for enhancement of power when the sparsity assumption is violated. 

In this paper, we introduce an ensemble approach and a binary expansion randomized ensemble test. The ensemble approach uses both the sum of squared symmetric statistics and the distance correlation test. It shows better power in linear and local settings while maintaining power for other dependence structures. Moreover, it can be easily generalized to an independence test for the multivariate setting, the binary expansion randomized ensemble test. By random projections, the BERET transforms the multivariate independence testing problem into a univariate testing problem. The BERET also maintains the clear interpretability of the maximum binary expansion testing.

Simulation studies suggest that the power of the BERET is advantageous compared with a range of competitors considered in many meaningful dependence structures. Investigation of three data examples shows that the BERET reveals hidden dependence structures from the data while maintaining a level of power similar to the best of the competing methods.

\vspace*{-10pt}

\section*{Acknowledgement}
The authors thank the helpful comments and suggestions from Xiao-Li Meng and Li Ma. This research was partially supported by the National Science Foundation and a grant from the National Cancer Institute. % DMS-1613112, NSF IIS-1633212, and NSF DMS-1916237

\counterwithin{theorem}{section}

\renewcommand\thesection{\Alph{section}}

\section*{Appendix: Proofs}

\newcommand\invisiblesection[1]{%
  \refstepcounter{section}%
  \addcontentsline{toc}{section}{\protect\numberline{\thesection}#1}%
  \sectionmark{#1}}
\invisiblesection{1}
\setcounter{section}{1}

\subsection{Proof of Theorem \ref{THM:uni_property}}

\begin{enumerate}[\rm(i)]
\item Suppose that $\mathcal{B}_d(X, Y)=0$. Then $E(\dot{A}_\mathbf{a}\dot{B}_\mathbf{b})^2\leq (2^d-1)^2\mathcal{B}_d(X, Y)=0$ for $\mathbf{a}\ne\mathbf{0}$ and $\mathbf{b}\ne\mathbf{0}$. If $E(\dot{A}_\mathbf{a}\dot{B}_\mathbf{b})=0$ for $\mathbf{a}\ne\mathbf{0}$ and $\mathbf{b}\ne\mathbf{0}$, then by definition, $\mathcal{B}_d(X, Y)=0$. Thus, by theorem 4.1 in \citet{zhang2019bet}, $\mathcal{B}_d(X, Y) = 0$ if and only if $U_d$ and $V_d$ are independent. \qed
\item Since $E(\dot{A}_\mathbf{a}\dot{B}_\mathbf{b})=2Pr(\dot{A}_\mathbf{a}\dot{B}_\mathbf{b}=1)-1$, we have $0\leq E(\dot{A}_\mathbf{a}\dot{B}_\mathbf{b})^2\leq 1$ and Therefore, $0\leq \sum_{\mathbf{ab}\in C_d}{E(\dot{A}_\mathbf{a}\dot{B}_\mathbf{b})^2}\leq (2^d-1)^2$. \qed
\item By the definition of $S_{(\mathbf{ab})}$, we obtain $0 \leq (S_{(\mathbf{ab})}/n)^2 \leq 1$. Since $\vert C_d\vert = (2^d-1)^2$, we have $0\leq \sum_{\mathbf{ab}\in C_d}{(S_{(\mathbf{ab})}/n)^2}\leq (2^d-1)^2$. \qed
\item By law of large numbers, $S_{\mathbf{(ab)}}/n$ converges almost surely to $E(\dot{A}_\mathbf{a}\dot{B}_\mathbf{b})$ for $\mathbf{a}\ne\mathbf{0}$ and $\mathbf{b}\ne\mathbf{0}$. Hence, the conclusion follows at once by the continuous mapping theorem.  \qed
\item Let $\mathbf{S}$ be a vector with entries the $S_{(\mathbf{ab})}$'s. Each $S_{(\mathbf{ab})}/n$ is a sample mean of $\dot{A}_\mathbf{a}\dot{B}_\mathbf{b}$ terms with mean 0 and variance 1. Since the $S_{(\mathbf{ab})}$'s are pairwise independent, by the central limit theorem, $\mathbf{S}/\sqrt{n}$ converges to $\mathcal{N}(0, \mathbf{I}_{(2^d-1)^2})$. By the continuous mapping theorem, each $(2^d-1)^2n\mathcal{B}_{n,d}[\{(X_i, Y_i)\}_{i=1}^n]$ is asymptotically $\chi^2$ with $(2^d - 1)^2$ degree of freedom. \qed

\end{enumerate}

\subsection{Proof of Theorem \ref{THM:uni_consistency}}

Suppose $TV(\mathbf{P}_{(U_d, V_d)}, \mathbf{P}_{0,d})\geq \delta$ for a fixed $0<\delta\leq 1/2$. Then, for some $(\mathbf{a}'\mathbf{b}')$, $E[\dot{A}_{\mathbf{a}'}\dot{B}_{\mathbf{b}'}]\geq \frac{2\sqrt{2d}\delta}{2^{d/2}}$ (see the proof of theorem 4.4 in \citet{zhang2019bet}). Therefore, we have $(2^d-1)^2\mathcal{B}_d(X, Y)\geq \frac{d\delta^2}{2^{d-3}}$. Since $\mathcal{B}_{n,d}[\{(X_i, Y_i)\}_{i=1}^n]\overset{a.s.}{\longrightarrow}\mathcal{B}_d(X, Y)$ by theorem \ref{THM:uni_property} (iv), $(2^d-1)^2n\mathcal{B}_{n,d}[\{(X_i, Y_i)\}_{i=1}^n]\rightarrow \infty~as~n\rightarrow \infty$. \qed

\subsection{Proof of Lemma \ref{LEM:Lemma1}}

Suppose $\mathbf{a}^T\mathbf{X}\perp \mathbf{b}^T\mathbf{Y}$ for all $\mathbf{a}\in \mathbb{R}^p, \mathbf{b}\in \mathbb{R}^q$ such that $\lVert \mathbf{a} \rVert=1$ and $\lVert \mathbf{b} \rVert=1$, we have
\begin{align*}
	E\big[\exp\{is(a_1X_1+&\cdots+a_pX_p)+it(b_1Y_1+\cdots+b_qY_q)\}\big] \\
	&=\phi_{\mathbf{a}^T\mathbf{X},{ \mathbf{b}^T\mathbf{Y}}}(s, t) \\
	&=\phi_{\mathbf{a}^T\mathbf{X}}(s)\phi_{\mathbf{b}^T\mathbf{Y}}(t) \\
	&=E\big[\exp\{is(a_1X_1+\cdots+a_pX_p)\}\big]E\big[\exp\{it(b_1Y_1+\cdots+b_qY_q)\}\big]
\end{align*}
for all $\mathbf{a}\in \mathbb{R}^p, \mathbf{b}\in \mathbb{R}^q$ such that $\lVert \mathbf{a} \rVert=1$ and $\lVert \mathbf{b} \rVert=1$. Now, consider the characteristic function of $\mathbf{X}$ and $\mathbf{Y}$. Then, by the above result, we obtain
\begin{eqnarray*}
	\phi_{\mathbf{X}, \mathbf{Y}}(\mathbf{s},\mathbf{t})&=&E\big\{\exp(i\mathbf{s}^T \mathbf{X}+i\mathbf{t}^T \mathbf{Y})\big\}\\
	&=& E\bigg[\exp\Big\{i\lVert \mathbf{s} \rVert\Big(\frac{s_1}{\lVert \mathbf{s} \rVert}X_1+\cdots+\frac{s_p}{\lVert \mathbf{s} \rVert}X_p\Big)+i\lVert \mathbf{t} \rVert\Big(\frac{t_1}{\lVert \mathbf{t} \rVert}Y_1+\cdots+\frac{t_q}{\lVert \mathbf{t} \rVert}Y_q\Big)\Big\}\bigg]\\
	&=& E\big\{\exp(i\mathbf{s}^T\mathbf{X})\big\}E\big\{\exp(i\mathbf{t}^T\mathbf{Y})\big\}\\
	&=& \phi_{X}(\mathbf{s})\phi_{Y}(\mathbf{t}).
\end{eqnarray*}
The opposite direction can also be easily shown.	\qed

\subsection{Proof of Theorem \ref{THM:multi_property}}

To prove theorem \ref{THM:multi_property}, we need the following lemmas.

\begin{lemma}
	Let $\mathbf{X}\in\mathbb{R}^p$ be a vector of $p$ continuous random variables and $\mathbf{s}_1, \mathbf{s}_2$ be two vectors in $\mathbb{R}^p$. Then $\lvert \mathbf{s}_1^T\mathbf{X}- \mathbf{s}_2^T\mathbf{X}\rvert\overset{P}{\longrightarrow}0$ as $\lVert \mathbf{s}_1 - \mathbf{s}_2\lVert\rightarrow 0$.
\label{LEM:lemma2}	
\end{lemma}
\noindent \textit{Proof.} Let $\epsilon>0$ and $\eta = \lVert \mathbf{s}_1 - \mathbf{s}_2\lVert$, we obtain
\begin{align*}
	Pr(\lvert \mathbf{s}_1^T\mathbf{X}- \mathbf{s}_2^T\mathbf{X}\rvert > \epsilon) &\leq Pr(\lVert \mathbf{s}_1 - \mathbf{s}_2\lVert\lVert \mathbf{X} \rVert>\epsilon)\\
	&\leq Pr(\eta\lVert \mathbf{X} \rVert>\epsilon)\\
	&\leq Pr(\lVert \mathbf{X} \rVert>\epsilon/\eta)\rightarrow 0~as~\eta\rightarrow 0.
\end{align*}
\qed

\begin{lemma}
	Let $U_i:=F_{\mathbf{s}_i^T\mathbf{X}}(\mathbf{s_i}^T\mathbf{X})$ for $i = 1, 2$. $\lvert U_1-U_2\rvert\overset{P}{\longrightarrow}0$ as $\lVert \mathbf{s}_1 - \mathbf{s}_2\lVert\rightarrow 0$.
\label{LEM:lemma3}
\end{lemma}

\noindent \textit{Proof.} \begin{align*}
	Pr(\lvert U_1-U_2\rvert>\epsilon)=&Pr(\lvert F_{\mathbf{s}_1^T\mathbf{X}}(\mathbf{s}_1^T\mathbf{X})-F_{\mathbf{s}_2^T\mathbf{X}}(\mathbf{s}_2^T\mathbf{X})\rvert>\epsilon)\\
	\leq& Pr(\lvert F_{\mathbf{s}_1^T\mathbf{X}}(\mathbf{s}_1^T\mathbf{X})-F_{\mathbf{s}_1^T\mathbf{X}}(\mathbf{s}_2^T\mathbf{X})\rvert>\epsilon/2) \\ &+Pr(\lvert F_{\mathbf{s}_1^T\mathbf{X}}(\mathbf{s}_2^T\mathbf{X})-F_{\mathbf{s}_2^T\mathbf{X}}(\mathbf{s}_2^T\mathbf{X})\rvert>\epsilon/2).
\end{align*}

Since $F_{\mathbf{s}_1^T\mathbf{X}}$ is uniformly continuous, there exist $\delta_1>0$ such that $\lvert F_{\mathbf{s}_1^T\mathbf{X}}(w_1)-F_{\mathbf{s}_1^T\mathbf{X}}(w_2)\rvert \leq \epsilon/2$ for all $w_1$ and $w_2$ with $|w_1-w_2| \leq \delta_1$. 
\begin{align*}
	Pr\{\lvert F_{\mathbf{s}_1^T\mathbf{X}}(\mathbf{s}_1^T\mathbf{X})-&F_{\mathbf{s}_1^T\mathbf{X}}(\mathbf{s}_2^T\mathbf{X})\rvert>\epsilon/2\} \\
	=&~Pr\{\lvert F_{\mathbf{s}_1^T\mathbf{X}}(\mathbf{s}_1^T\mathbf{X})-F_{\mathbf{s}_1^T\mathbf{X}}(\mathbf{s}_2^T\mathbf{X})\rvert>\epsilon/2, \lvert \mathbf{s}_1^T\mathbf{X}- \mathbf{s}_2^T\mathbf{X}\rvert \leq \delta_1\} \\
	&+Pr\{\lvert F_{\mathbf{s}_1^T\mathbf{X}}(\mathbf{s}_1^T\mathbf{X})-F_{\mathbf{s}_1^T\mathbf{X}}(\mathbf{s}_2^T\mathbf{X})\rvert>\epsilon/2, \lvert \mathbf{s}_1^T\mathbf{X}- \mathbf{s}_2^T\mathbf{X}\rvert > \delta_1\} \\
	\leq&~Pr(\lvert \mathbf{s}_1^T\mathbf{X}- \mathbf{s}_2^T\mathbf{X}\rvert > \delta_1).
\end{align*}

Now fix $w\in \mathbb{R}$.
\begin{align*}
	\lvert F_{\mathbf{s}_1^T\mathbf{X}}(w)-F_{\mathbf{s}_2^T\mathbf{X}}(w)\rvert \leq &~\lvert Pr(\mathbf{s}_1^T\mathbf{X}\leq w, \mathbf{s}_2^T\mathbf{X}\leq w)+Pr(\mathbf{s}_1^T\mathbf{X}\leq w, \mathbf{s}_2^T\mathbf{X}> w)\\
	&-Pr(\mathbf{s}_1^T\mathbf{X}\leq w, \mathbf{s}_2^T\mathbf{X}\leq w)-Pr(\mathbf{s}_1^T\mathbf{X}> w, \mathbf{s}_2^T\mathbf{X}\leq w)\rvert\\
	\leq&~Pr(\mathbf{s}_1^T\mathbf{X}\leq w, \mathbf{s}_2^T\mathbf{X}> w) + Pr(\mathbf{s}_1^T\mathbf{X}> w, \mathbf{s}_2^T\mathbf{X}\leq w) \\
	\leq&~Pr(\mathbf{s}_1^T\mathbf{X}\leq w, \mathbf{s}_2^T\mathbf{X}> w, \lvert \mathbf{s}_1^T\mathbf{X}- \mathbf{s}_2^T\mathbf{X} \rvert \leq\delta_2)\\
	&+Pr(\mathbf{s}_1^T\mathbf{X}\leq w, \mathbf{s}_2^T\mathbf{X}> w, \lvert \mathbf{s}_1^T\mathbf{X}- \mathbf{s}_2^T\mathbf{X} \rvert>\delta_2)\\
	&+Pr(\mathbf{s}_1^T\mathbf{X}> w, \mathbf{s}_2^T\mathbf{X}\leq w, \lvert \mathbf{s}_1^T\mathbf{X}- \mathbf{s}_2^T\mathbf{X} \rvert\leq\delta_2)\\
	&+Pr(\mathbf{s}_1^T\mathbf{X}> w, \mathbf{s}_2^T\mathbf{X}\leq w, \lvert \mathbf{s}_1^T\mathbf{X}- \mathbf{s}_2^T\mathbf{X} \rvert>\delta_2)\\
	\leq&~Pr(\mathbf{s}_2^T\mathbf{X}-\delta_2\leq w\leq\mathbf{s}_2^T\mathbf{X}+\delta_2)+Pr(\mathbf{s}_1^T\mathbf{X}-\delta_2\leq w\leq\mathbf{s}_1^T\mathbf{X}+\delta_2)\\
	&+Pr(\lvert \mathbf{s}_1^T\mathbf{X}- \mathbf{s}_2^T\mathbf{X} \rvert>\delta_2)	\\
	\leq&~Pr(w-\delta_2\leq\mathbf{s}_2^T\mathbf{X}\leq w+\delta_2)+Pr(w-\delta_2\leq\mathbf{s}_1^T\mathbf{X}\leq w+\delta_2)\\
	&+Pr(\lvert \mathbf{s}_1^T\mathbf{X}- \mathbf{s}_2^T\mathbf{X} \rvert>\delta_2).
\end{align*}

Since $F_{\mathbf{s}_1^T\mathbf{X}}$ and $F_{\mathbf{s}_2^T\mathbf{X}}$ are uniformly continuous and $\lvert \mathbf{s}_1^T\mathbf{X}- \mathbf{s}_2^T\mathbf{X}\rvert\overset{P}{\longrightarrow}0$, we can find $\delta_2$ such that the right hand side of the last inequality is less than $\epsilon/2$ for all $w\in \mathbb{R}$. Thus, we obtain
\begin{align*}
	Pr(\lvert F_{\mathbf{s}_1^T\mathbf{X}}(\mathbf{s}_2^T\mathbf{X})-&F_{\mathbf{s}_2^T\mathbf{X}}(\mathbf{s}_2^T\mathbf{X})\rvert>\epsilon/2) \\
	=&~Pr\{\lvert F_{\mathbf{s}_1^T\mathbf{X}}(\mathbf{s}_2^T\mathbf{X})-F_{\mathbf{s}_2^T\mathbf{X}}(\mathbf{s}_2^T\mathbf{X})\rvert>\epsilon/2, \lvert \mathbf{s}_1^T\mathbf{X}- \mathbf{s}_2^T\mathbf{X}\rvert \leq \delta_2\} \\
	&+Pr\{\lvert F_{\mathbf{s}_1^T\mathbf{X}}(\mathbf{s}_2^T\mathbf{X})-F_{\mathbf{s}_2^T\mathbf{X}}(\mathbf{s}_2^T\mathbf{X})\rvert>\epsilon/2, \lvert \mathbf{s}_1^T\mathbf{X}- \mathbf{s}_2^T\mathbf{X}\rvert > \delta_2\} \\
	\leq&~Pr(\lvert \mathbf{s}_1^T\mathbf{X}- \mathbf{s}_2^T\mathbf{X}\rvert > \delta_2).
\end{align*}
Therefore we have
\begin{align*}
	Pr(\lvert U_1-U_2\rvert>\epsilon)\leq&~Pr(\lvert F_{\mathbf{s}_1^T\mathbf{X}}(\mathbf{s}_1^T\mathbf{X})-F_{\mathbf{s}_1^T\mathbf{X}}(\mathbf{s}_2^T\mathbf{X})\rvert>\epsilon/2) \\ &+Pr(\lvert F_{\mathbf{s}_1^T\mathbf{X}}(\mathbf{s}_2^T\mathbf{X})-F_{\mathbf{s}_2^T\mathbf{X}}(\mathbf{s}_2^T\mathbf{X})\rvert>\epsilon/2)\\
	\leq&~Pr(\lvert \mathbf{s}_1^T\mathbf{X}- \mathbf{s}_2^T\mathbf{X}\rvert > \delta_1)+Pr(\lvert \mathbf{s}_1^T\mathbf{X}- \mathbf{s}_2^T\mathbf{X}\rvert > \delta_2).
\end{align*}

By Lemma \ref{LEM:lemma2}, the right hand side of the last inequality vanishes as $\lVert \mathbf{s}_1 - \mathbf{s}_2\lVert\rightarrow 0$.
\qed

\begin{lemma}
	For any given $\epsilon>0$, an integer $d>0$ and two intervals $I_{i}=(\frac{i-1}{2^d}, \frac{i}{2^d}], I_{j}=(\frac{j-1}{2^d}, \frac{j}{2^d}]$ with $1\leq i, j, \leq 2^d$, there exists $\eta>0$ such that $\lvert E\{(I(U_1\in I_{i}, V_1\in I_{j})\} - E\{I(U_2\in I_{i}, V_2\in I_{j})\}\rvert <\epsilon$ whenever $\lVert \mathbf{s}_1-\mathbf{s}_2\rVert<\eta$ and $\lVert \mathbf{t}_1-\mathbf{t}_2\rVert<\eta$.
\label{LEM:lemma4}
\end{lemma}
 
\noindent \textit{Proof.} For any $0\leq a<b\leq 1$, we have
\begin{align*}
	E\{I(U_1\in (a, b], U_2\notin (a, b])\} =&~Pr(U_1\in (a, b], U_2\notin (a, b]) \\
	=&~Pr(U_1\in (a, b], U_2\notin (a, b], \lvert U_1-U_2 \rvert < \delta/2)\\
	&+ Pr(U_1\in (a, b], U_2\notin (a, b], \lvert U_1-U_2 \rvert \geq \delta/2) \\
	\leq&~\delta + Pr(\lvert U_1-U_2 \rvert \geq \delta/2)
\end{align*}
Since $\delta$ is arbitrary, by Lemma \ref{LEM:lemma3}, we can find $\eta > 0$ such that $E\{I(U_1\in (a, b], U_2\notin (a, b])\}\leq \epsilon/4$ whenever $\lVert \mathbf{s}_1-\mathbf{s}_2\rVert<\eta$. In the same manner, we can bound $E\{I(V_1\in (a', b'], V_2\notin (a', b'])\}$ for any $0\leq a'<b'\leq 1$. Now for any $I_{i}=(\frac{i-1}{2^d}, \frac{i}{2^d}], I_{j}=(\frac{j-1}{2^d}, \frac{j}{2^d}]$ with $1\leq i, j, \leq 2^d$, we obtain
\begin{align*}
	\lvert E\{(I(U_1\in I_{i}, V_1\in I_{j})\} &- E\{I(U_2\in I_{i}, V_2\in I_{j})\}\rvert \\
	=&~\lvert E\{(I(U_1\in I_{i}, V_1\in I_{j})-I(U_2\in I_{i}, V_2\in I_{j})\}\rvert \\
	\leq&~\lvert E\{(I(U_1\in I_{i}, U_2\notin I_{i})\}\rvert + \lvert E\{I(U_2\in I_{i}, V_2\notin I_{j})\}\rvert \\
	&~+\lvert E\{(I(U_1\notin I_{i}, U_2\in I_{i})\}\rvert + \lvert E\{I(U_2\notin I_{i}, V_2\in I_{j})\}\rvert
\end{align*}
By the above results, we can find $\eta_1, \eta_2, \eta_3$ and $\eta_4$ that bound each term by $\epsilon/4$ in the last inequality. By letting $\eta = min\{\eta_1, \eta_2, \eta_3, \eta_4\}$, we can have the desired result.
\qed

\begin{lemma}
	If $\mathcal{B}_d(\mathbf{s_0}^T\mathbf{X}, \mathbf{t_0}^T\mathbf{Y})>c$, then there exist $\delta > 0$ such that $\mathcal{B}_d(\mathbf{s}^T\mathbf{X}, \mathbf{t}^T\mathbf{Y})>c/2$ for all $(\mathbf{s}, \mathbf{t})\in B_\delta(\mathbf{s}_0, \mathbf{t}_0):=\{(\mathbf{s}, \mathbf{t})\mid \lVert\mathbf{s}-\mathbf{s}_0\rVert^2+\lVert\mathbf{t}-\mathbf{t}_0\rVert^2\leq\delta~and~\lVert \mathbf{s} \rVert=\lVert \mathbf{t} \rVert=1\}$.	
\label{LEM:lemma5}
\end{lemma}
\noindent \textit{Proof.} Fix $d>0$. Each cross interaction can be expressed as a sum of $2^{2d}$ signed indicator variables, that is, 
\begin{align*}
	\dot{A}_{\mathbf{a}}(\mathbf{s}^T\mathbf{X})\dot{B}_{\mathbf{b}}(\mathbf{t}^T\mathbf{Y})=\sum_{i=1}^{2^d}{\sum_{j=1}^{2^d}{(-1)^{h(i,j)}I\Big\{U\in \Big(\frac{i-1}{2^d}, \frac{i}{2^d}\Big], V\in \Big(\frac{j-1}{2^d}, \frac{j}{2^d}\Big]\Big\}}},
\end{align*}
where $U_{\mathbf{s}}=F_{\mathbf{s}^T\mathbf{X}}(\mathbf{s}^T\mathbf{X})$, $V_{\mathbf{t}}=F_{\mathbf{t}^T\mathbf{Y}}(\mathbf{t}^T\mathbf{Y})$ and $h(i, j)\in \{0, 1\}$. Let $\epsilon = \frac{c}{2^{2d+2}}$. By Lemma \ref{LEM:lemma4}, we can find $\eta_{i,j}$ for each $(i, j)$ such that $\lvert E\{(I(U_{\mathbf{s}}\in I_{i}, V_{\mathbf{t}}\in I_{j})\} - E\{I(U_{\mathbf{s}_0}\in I_{i}, V_{\mathbf{t}_0}\in I_{j})\}\rvert <\epsilon$. By letting $\delta_{\mathbf{ab}}=\min_{1\leq i,j\leq 2^d}{\eta_{i,j}^2}$, we can obtain
\begin{align*}
	\lvert E\{\dot{A}_{\mathbf{a}}(\mathbf{s}^T\mathbf{X})\dot{B}_{\mathbf{b}}(\mathbf{t}^T\mathbf{X})\}-E\{\dot{A}_{\mathbf{a}}(\mathbf{s_0}^T\mathbf{X})\dot{B}_{\mathbf{b}}(\mathbf{t_0}^T\mathbf{X})\}\rvert\leq\frac{c}{4},
\end{align*}
for all $(\mathbf{s}, \mathbf{t})\in B_{\delta_{\mathbf{ab}}}(\mathbf{s}_0, \mathbf{t}_0)=\{(\mathbf{s}, \mathbf{t})\lvert \lVert\mathbf{s}-\mathbf{s}_0\rVert^2+\lVert\mathbf{t}-\mathbf{t}_0\rVert^2\leq\delta_{\mathbf{ab}}~and~\lVert \mathbf{s} \rVert=\lVert \mathbf{t} \rVert=1\}$.
By definition, population measure of dependence is sum of squared expectations of cross interactions, that is,
\begin{align*}
	\mathcal{B}_d(\mathbf{s}^T\mathbf{X}, \mathbf{t}^T\mathbf{Y})=\frac{1}{(2^d-1)^2}\sum_{ab\in C_d}{E[\dot{A}_{\mathbf{a}}\dot{B}_{\mathbf{b}}]^2}.
\end{align*}
By letting $\delta=\min_{\mathbf{ab}\in C_d}{\delta_{\mathbf{ab}}}$, we obtain
\begin{align*}
	\lvert \mathcal{B}_d(\mathbf{s}^T\mathbf{X}, \mathbf{t}^T\mathbf{Y})&-\mathcal{B}_d(\mathbf{s_0}^T\mathbf{X}, \mathbf{t_0}^T\mathbf{Y})\}\rvert\\
	&=\frac{1}{(2^d-1)^2}\sum_{ab\in C_d}{\lvert E\{\dot{A}_{\mathbf{a}}(\mathbf{s}^T\mathbf{X})\dot{B}_{\mathbf{b}}(\mathbf{t}^T\mathbf{X})\}^2-E\{\dot{A}_{\mathbf{a}}(\mathbf{s_0}^T\mathbf{X})\dot{B}_{\mathbf{b}}(\mathbf{t_0}^T\mathbf{X})\}^2\rvert}\\
	&\leq\frac{1}{(2^d-1)^2}\sum_{ab\in C_d}{2\lvert E\{\dot{A}_{\mathbf{a}}(\mathbf{s}^T\mathbf{X})\dot{B}_{\mathbf{b}}(\mathbf{t}^T\mathbf{X})\}-E\{\dot{A}_{\mathbf{a}}(\mathbf{s_0}^T\mathbf{X})\dot{B}_{\mathbf{b}}(\mathbf{t_0}^T\mathbf{X})\}\rvert}\\
	&\leq\frac{c}{2}
\end{align*}
for all $(\mathbf{s}, \mathbf{t})\in B_\delta(\mathbf{s}_0, \mathbf{t}_0)=\{(\mathbf{s}, \mathbf{t})\lvert \lVert\mathbf{s}-\mathbf{s}_0\rVert^2+\lVert\mathbf{t}-\mathbf{t}_0\rVert^2\leq\delta~and~\lVert \mathbf{s} \rVert=\lVert \mathbf{t} \rVert=1\}$.
Hence, we can have desired result.
\qed \\

The proof of theorem \ref{THM:multi_property} is given in the following list.

\begin{enumerate}[\rm(i)]
	\item Suppose $U_d^{\mathbf{s}}$ and $V_d^{\mathbf{t}}$ are independent for all $\mathbf{s}\in S_p, \mathbf{t}\in S_q$. By theorem \ref{THM:uni_property} (i), $\mathcal{B}_d(\mathbf{s}^T\mathbf{X}, \mathbf{t}^T\mathbf{Y})=0$ for all $\mathbf{s}\in S_p, \mathbf{t}\in S_q$. Therefore, we obtain $\mathcal{B}_d(\mathbf{X}, \mathbf{Y})=0$. Now suppose that there exists $(\mathbf{s}_0, \mathbf{t}_0)$ such that $U_d^{\mathbf{s_0}}$ and $V_d^{\mathbf{t_0}}$ are not independent. Let $c=\mathcal{B}_d(\mathbf{s_0}^T\mathbf{X}, \mathbf{t_0}^T\mathbf{Y})>0$. Then, by Lemma \ref{LEM:lemma5}, there exist $\delta>0$ such that $\mathcal{B}_d(\mathbf{s}^T\mathbf{X}, \mathbf{t}^T\mathbf{Y})>c/2$ for all $(\mathbf{s}, \mathbf{t})\in B_\delta(\mathbf{s}_0, \mathbf{t}_0)$. Therefore, we obtain
	\begin{align*}
		\mathcal{B}_d(\mathbf{X}, \mathbf{Y}) &= \frac{1}{c_pc_q}\int_{S_q}{\int_{S_p}{\mathcal{B}_d(\mathbf{s}^T\mathbf{X}, \mathbf{t}^T\mathbf{Y})d\mathbf{s}d\mathbf{t}}}\\&>\frac{c}{2c_p c_q}\mathcal{L}\{B_\delta(\mathbf{s}_0, \mathbf{t}_0)\}\\&>0.
	\end{align*}
	\qed
	\item The proof is immediate and is omitted.
	\item The proof is immediate and is omitted.
	\item Since $\mathcal{B}_{n, d}(\{(\mathbf{s}^T\mathbf{X}_i, \mathbf{t}^T\mathbf{Y}_i)\}_{i=1}^n)\overset{a.s.}{\longrightarrow}\mathcal{B}_d(\mathbf{s}^T\mathbf{X}, \mathbf{t}^T\mathbf{Y})$ as $n\rightarrow\infty$ by theorem \ref{THM:uni_property} (v), $\mathcal{B}_{n, d}(\{(\mathbf{X}_i, \mathbf{Y}_i)\}_{i=1}^n)$ converges to $\mathcal{B}_d(\mathbf{X}, \mathbf{Y})$ by dominated convergence theorem.	\qed
\end{enumerate}

\subsection{Proof of Theorem \ref{THM:multi_convergence}}

To prove theorem \ref{THM:multi_convergence}, we need following lemmas and corollary.

\begin{lemma}
$\mathcal{F}_1=\{I(\mathbf{s}^T\mathbf{X}\leq u):\mathbf{s}\in S_p, u\in \mathbb{R}\}$ is a Donsker class.
\label{THM:Lemma6}
\end{lemma}
\noindent \textit{Proof.} This folllows directly from lemmas 8.12, 9.6, and 9.9 of \citet{kosorok2007introduction}. \qed

\begin{lemma}
Let $\mathcal{F}_{*}$ be a class of all cadlag cumulative distribution functions $F$ with $\lim_{x\rightarrow \infty}{F(x)=1}$ and $ \lim_{x\rightarrow-\infty}{F(x)}=0$. Then $\mathcal{F}_2=\{I\{F(\mathbf{s}^T\mathbf{X})\leq u\}:F \in\mathcal{F}_*, \mathbf{s}\in S_p, u\in \mathbb{R}\}$ is a Donsker class.
\label{THM:Lemma7}
\end{lemma}

\noindent \textit{Proof.} Note that for any $F\in \mathcal{F}_*$, $I\{F(\mathbf{s}^T\mathbf{X})\leq u\}=I\{\mathbf{s}^T\mathbf{X}\leq F^{-1}(u)\}$ a.s. for all $\mathbf{s}\in S_p$ and all $u\in\mathbb{R}$. Thus $\mathcal{F}_2=\{I\{\mathbf{s}^T\mathbf{X}\leq F^{-1}(u)\}:\mathbf{s}\in S_p, u\in\mathbb{R}, F\in \mathcal{F}_*\}\subset\{I\{\mathbf{s}^T\mathbf{X}\leq z\}:\mathbf{s}\in S_p, z\in\mathbb{R}\}=\mathcal{F}_1$. Therefore, by Lemma \ref{THM:Lemma6}, $\mathcal{F}_2$ is a Donsker class. \qed

\begin{corollary}
The following sets of functions are Donsker classes:
\begin{enumerate}[\rm(i)]
	\item $\mathcal{F}_3=\{I\{u_1 < F(\mathbf{s}^T\mathbf{X})\leq u_2\}:F \in\mathcal{F}_*, \mathbf{s}\in S_p, u_1, u_2\in \mathbb{R}\}$
	\item $\mathcal{F}_4=\{I\{F(\mathbf{s}^T\mathbf{X})\in A\}:F \in\mathcal{F}_*, \mathbf{s}\in S_p, A\in \mathcal{A}\}$
	\item $\mathcal{F}_5=\{I\{F(\mathbf{s}^T\mathbf{X})\in A, F(\mathbf{t}^T\mathbf{Y})\in B\}:F, G \in\mathcal{F}_*, s\in S_p, \mathbf{t}\in S_q, A\in \mathcal{A}, B\in \mathcal{B}\}$
\end{enumerate}
\label{THM:corollary1}
\end{corollary}

\noindent \textit{Proof.} 
\begin{enumerate}[\rm(i)]
	\item It follows from Lemma \ref{THM:Lemma7} and the fact that products of bounded Donsker classes are Donsker classes.\qed
	\item It follows from (i) since $\mathcal{A}$ is collection of unions and intersections of intervals and reapplication of the preservation of the Donsker property under products of bounded Donsker classes.\qed
	\item It follows from (ii) and reapplication of the preservation of the Donsker property under products of bounded Donsker classes.\qed
\end{enumerate}

By Corollary \ref{THM:corollary1} and Glivenko-Cantelli theorem, we obtain
 \begin{align}
 	\sup_{\mathbf{s}\in S_p,\mathbf{t}\in S_q}\biggl\lvert n^{-1}\sum_{i=1}^n{I\{\widehat{F}_{n,\mathbf{s}}(\mathbf{s}^T\mathbf{X}_i)}&\in A, \widehat{G}_{n,\mathbf{t}}(\mathbf{t}^T\mathbf{Y}_i)\in B\} \nonumber \\ &-\mathbb{E}_{X,Y}\big[I\{\widehat{F}_{n,\mathbf{s}}(\mathbf{s}^T\mathbf{X})\in A, \widehat{G}_{n,\mathbf{t}}(\mathbf{t}^T\mathbf{Y})\in B\big]\biggr\rvert\overset{a.s.}{\longrightarrow}0
 \label{conv5}
 \end{align}
where $\widehat{F}_{n,\mathbf{s}}$ is the empirical cdf of $\mathbf{s}^T\mathbf{X}$ and $\widehat{G}_{n,\mathbf{t}}$ is that of $\mathbf{t}^T\mathbf{Y}$.

Note that $\widehat{F}_{n,\mathbf{s}}(u)=n^{-1}\sum_{i=1}^n{I(\mathbf{s}^T\mathbf{X}_i\leq u})$ and $\widehat{G}_{n,\mathbf{t}}(v)=n^{-1}\sum_{i=1}^n{I(\mathbf{t}^T\mathbf{Y}_i\leq v})$, recycling previous arguments, we have
\begin{align*}
	\sup_{\mathbf{s}\in S_p, u\in \mathbb{R}}{\lvert \widehat{F}_{n,\mathbf{s}}(u)-F_{0,s}(u)\rvert} + \sup_{\mathbf{t}\in S_q, v\in \mathbb{R}}{\lvert \widehat{G}_{n,\mathbf{t}}(v)-G_{0,t}(v)\rvert}\overset{a.s.}{\longrightarrow}0,
\end{align*}
where $F_{0, s}(u)=Pr(\mathbf{s}^T\mathbf{X}\leq u)$ and $G_{0, t}(v)=Pr(\mathbf{t}^T\mathbf{Y}\leq v)$.

Assume $X$ and $Y$ are continuous. And for any $\epsilon>0$, $A\subset \mathbb{R}$, define $A^{\epsilon}=\{x\in\mathbb{R}:\inf_{y\in A}{\lVert x- y\rVert<\epsilon}\}$ and $A_{\epsilon}=((A^c)^\epsilon)^c$. Note that by continuity of $X$ and $Y$,
\begin{align}
	\lim_{\epsilon\rightarrow0}{\sup_{\mathbf{s}\in S_p}{\lvert Pr(\mathbf{s}^T\mathbf{X}\in A^{\epsilon})-Pr(\mathbf{s}^T\mathbf{X}\in A)\rvert}}\rightarrow0,
\label{conv1}
\end{align}
similarly, for $A_{\epsilon}$ and also for $Y$ and $\mathbf{t}\in S_q$.

Now, fix $\epsilon>0$, by (\ref{conv5}), we have that for any $\mathbf{s}\in S_p, \mathbf{t}\in S_q$, 
\begin{align*}
	n^{-1}\sum_{i=1}^n{I\{\widehat{F}_{n,\mathbf{s}}(\mathbf{s}^T\mathbf{X}_i)\in A, \widehat{G}_{n,\mathbf{t}}(\mathbf{t}^T\mathbf{Y}_i)\in B\}}\leq&~n^{-1}\sum_{i=1}^n{I\{F_{0,s}(\mathbf{s}^T\mathbf{X})\in A^{\epsilon}, G_{0,t}(\mathbf{t}^T\mathbf{Y})\in B^{\epsilon}\}}\\&(\mathrm{for~all~n~large~enough~almost~surely})\\&\overset{a.s.}{\longrightarrow}\mathbb{E}\big[I\{F_{0,s}(\mathbf{s}^T\mathbf{X})\in A^{\epsilon}, G_{0,t}(\mathbf{t}^T\mathbf{Y})\in B^{\epsilon}\}\big]
\end{align*}
Similarly,
\begin{align*}
	n^{-1}\sum_{i=1}^n{I\{\widehat{F}_{n,\mathbf{s}}(\mathbf{s}^T\mathbf{X}_i)\in A, \widehat{G}_{n,\mathbf{t}}(\mathbf{t}^T\mathbf{Y}_i)\in B\}}\geq&~n^{-1}\sum_{i=1}^n{I\{F_{0,s}(\mathbf{s}^T\mathbf{X})\in A_{\epsilon}, G_{0,t}(\mathbf{t}^T\mathbf{Y})\in B_{\epsilon}\}}\\&\overset{a.s.}{\longrightarrow}\mathbb{E}\big[I\{F_{0,s}(\mathbf{s}^T\mathbf{X})\in A_{\epsilon}, G_{0,t}(\mathbf{t}^T\mathbf{Y})\in B_{\epsilon}\}\big]
\end{align*}
By (\ref{conv1}) and the fact that $\epsilon$ is arbitrary, uniformly for any $\mathbf{s}\in S_p, \mathbf{t}\in S_q$, we have
 \begin{align}
 	\sup_{\mathbf{s}\in S_p,\mathbf{t}\in S_q}\biggl\lvert n^{-1}\sum_{i=1}^n{I\bigl\{\widehat{F}_{n,\mathbf{s}}(\mathbf{s}^T\mathbf{X}_i)}&{\in A, \widehat{G}_{n,\mathbf{t}}(\mathbf{t}^T\mathbf{Y}_i)\in B\bigr\}}\nonumber\\&{-\mathbb{E}_{X,Y}\big[I\{F_{0,s}(\mathbf{s}^T\mathbf{X})\in A, G_{0,t}(\mathbf{t}^T\mathbf{Y})\in B\big\}]}\biggr\rvert\overset{a.s.}{\longrightarrow}0
\label{conv2}
\end{align}
 
Let $\mathcal{B}_{n,d}(\mathbf{s}, \mathbf{t})=\mathcal{B}_{n,d}(\{(\mathbf{s}^T\mathbf{X}, \mathbf{t}^T\mathbf{Y})\}_{i=1}^n)$ and $\mathcal{B}_{0,d}(\mathbf{s}, \mathbf{t})=\mathcal{B}_{d}(\mathbf{s}^T\mathbf{X}, \mathbf{t}^T\mathbf{Y})$. Now (\ref{conv2}) implies
\begin{align}
\sup_{\mathbf{s}\in S_p, \mathbf{t}\in S_q}{\lvert\mathcal{B}_{n,d}(\mathbf{s},\mathbf{t})-\mathcal{B}_{0,d}(\mathbf{s},\mathbf{t})\rvert}\overset{a.s.}{\longrightarrow}0.
\label{conv3}	
\end{align}
Let $(\mathbf{S}_1, \mathbf{T}_1),\dots, (\mathbf{S}_m, \mathbf{T}_m)$ be iid uniformly distributed on $S_p \times S_q$. Then by continuity of $\mathcal{B}_{0,d}$ in $\mathbf{s}$ and $\mathbf{t}$, it is easy to show
\begin{align}
m^{-1}\sum_{j=1}^m{\mathcal{B}_{0,d}(\mathbf{S}_j, \mathbf{T}_j)\overset{a.s.}{\longrightarrow}\frac{1}{c_pc_q}\int_{S_q}\int_{S_p}{{\mathcal{B}_{0,d}(\mathbf{s},\mathbf{t})dsdt}}}.
\label{conv4}
\end{align}
 
For $\epsilon>0$, by (\ref{conv3}), there exists $n_0<\infty$ such that $Pr(\sup_{n\geq n_0}{\lvert \mathcal{B}_{n,d}(\mathbf{s},\mathbf{t})-\mathcal{B}_{0,d}(\mathbf{s},\mathbf{t}) \rvert}>\epsilon)<\epsilon$ and by (\ref{conv4}), there exists $m_0<\infty$ such that
\begin{align*}
Pr\Bigg(\sup_{m\geq m_0}{\bigg\lvert m^{-1}\sum_{j=1}^{m}{\mathcal{B}_{0,d}(\mathbf{S}_j, \mathbf{T}_j)}-\frac{1}{c_pc_q}\int_{S_q}\int_{S_p}{{\mathcal{B}_{0,d}(\mathbf{s},\mathbf{t})dsdt}} \bigg\rvert}>\epsilon\Bigg)<\epsilon	
\end{align*}

Thus,
\begin{align*}
	\bigg\lvert m^{-1}\sum_{j=1}^m{\mathcal{B}_{n,d}(\mathbf{S}_j, \mathbf{T}_j)}-T_0\bigg\rvert \leq &\bigg\lvert m^{-1}\sum_{j=1}^m{\{\mathcal{B}_{n,d}(\mathbf{S}_j, \mathbf{T}_j)-\mathcal{B}_{0,d}(\mathbf{S}_j, \mathbf{T}_j)\}}\bigg\rvert\\ 
	&+ \bigg\lvert m^{-1}\sum_{j=1}^m{\mathcal{B}_{0,d}(\mathbf{S}_j, \mathbf{T}_j)}-T_0\bigg\rvert\\
	\leq & \sup_{\mathbf{s}\in S_p, \mathbf{t}\in S_q}{\bigg\lvert \mathcal{B}_{n,d}(\mathbf{s}, \mathbf{t})-\mathcal{B}_{0,d}(\mathbf{s}, \mathbf{t})\bigg\rvert}\\
	&+ \bigg\lvert m^{-1}\sum_{j=1}^m{\mathcal{B}_{0,d}(\mathbf{S}_j, \mathbf{T}_j)}-T_0\bigg\rvert \overset{a.s.}{\longrightarrow}0~as~n,m\rightarrow \infty.
\end{align*}
\qed

\subsection{Proof of Theorem \ref{THM:multi_consistency}}

\begin{enumerate}[\rm(i)]
	\item Suppose $TV(\mathbf{P}_{(U_d^{\mathbf{s}}, V_d^{\mathbf{t}})},\mathbf{P}_{0,d})\geq \delta$ for a fixed $0<\delta\leq 1/2$ and for some $\mathbf{s}\in S_p, \mathbf{t}\in S_q$. Then, there exist $(\mathbf{a}'\mathbf{b}')$ and $(\mathbf{s}_0, \mathbf{t}_0)$ such that $E(\dot{A}_{\mathbf{a}'}\dot{B}_{\mathbf{b}'})\geq \frac{2\sqrt{2d}\delta}{2^{d/2}}$ (see the proof of theorem 4.4 in \citet{zhang2019bet}). Therefore, we have $(2^d-1)^2\mathcal{B}_d(\mathbf{s}_0^T\mathbf{X}, \mathbf{t}_0^T\mathbf{Y})\geq \frac{d\delta^2}{2^{d-3}}$. By lemma \ref{LEM:lemma5}, there exist $\epsilon>0$ such that $(2^d-1)^2\mathcal{B}_d(\mathbf{s}^T\mathbf{X}, \mathbf{t}^T\mathbf{Y})\geq \frac{d\delta^2}{2^{d-2}}$ for all $(\mathbf{s}, \mathbf{t})\in B_\epsilon(\mathbf{s}_0, \mathbf{t}_0)$.	 Hence $(2^d-1)^2\mathcal{B}_d(\mathbf{X}, \mathbf{Y})=c$ for some constant $c>0$ for fixed $\delta$. Then, by theorem \ref{THM:multi_convergence}, $(2^d-1)^2\widehat{\mathcal{B}}^m_{n, d}[\{(\mathbf{X}_i, \mathbf{Y}_i)\}_{i=1}^n]$ converges to $(2^d-1)^2\mathcal{B}_d(\mathbf{X}, \mathbf{Y})=c>0$. Therefore, $n(2^d-1)^2\mathcal{B}_{n, d}[\{(\mathbf{X}_i, \mathbf{Y}_i)\}_{i=1}^n]\rightarrow\infty$ as $m, n\rightarrow\infty$ \qed
	\item Recall that $\zeta^{m}_{n,d}$ is nonnegative and symmetric in the order in which the observations for $i=1,\ldots,n$ appear. Moreover, under $H_0$ ($\mathbf{X}$ and $\mathbf{Y}$ are independent), $\zeta^{m}_{n,d} \rightarrow 0$ almost surely, while under $H_1$, $\zeta^{m}_{n,d} \rightarrow c>0$ almost surely. Now let $\Sigma_n$ be the $\sigma$-field generated by the information contained in $X_1,\ldots,X_n$ and $Y_1,\ldots,Y_n$ minus the ordering of the data and the pairing. Now let $\zeta^{m\ast}_{n,d}$ be a permuted version of $\zeta^{m}_{n,d}$ (the pairing is randomly permuted).

Assume the null hypothesis $H_0$. It is easy to verify that, conditional on $\Sigma_{n}$, both $\zeta^{m}_{n,d}$ and $\zeta^{m\ast}_{n,d}$ have the identical distribution, and thus the permutation based p-value of $\zeta^{m}_{n,d}$ will be uniformly distributed over $\{1/(n_p + 1),\ldots,1\}$, where $n_p$ is the number of permutation used (assuming that there are no duplicates). Even if we just use a randomly selected subset of the possible permutations, the uniformity will be approximately true. In addition, for every $\epsilon>0$, $Pr(\zeta^{m\ast}_{n,d}>\epsilon) = E\{Pr(\zeta^{m\ast}_{n,d}>\epsilon|\Sigma_n)\}=E\{Pr(\zeta^{m}_{n,d}>\epsilon|\Sigma_{n})\}=Pr(\zeta^{m}_{n,d}>\epsilon)\rightarrow 0$ by our previously established properties of $\zeta^{m}_{n,d}$ under $H_0$.

Now consider the alternative hypothesis $H_1$. In this case $\zeta^{m}_{n,d}\rightarrow c>0$ almost surely. Let $\zeta^{m0}_{n,d}$ be the test statistics based on data with $\mathbf{X}$ and $\mathbf{Y}$ indepedent but with the marginal distributions of $\mathbf{X}$ and $\mathbf{Y}$ being the same as the corresponding marginals under $H_1$. Let $\zeta^{m\ast}_{n,d}$ be the permutation version of the test statistics under $H_1$. Now it is easy to verify, recycling previous arguments, that $\zeta^{m\ast}_{n,d}$ and $\zeta^{m0}_{n,d}$ have identical distributions conditional on $\Sigma_n$. Thus we have once again, for every $\epsilon>0$, that $Pr(\zeta^{m\ast}_{n,d}>\epsilon)\rightarrow 0$ in probability. However, since $\zeta^{m}_{n,d}\rightarrow c>0$, we have that the permutation p-value of $\zeta^{m}_{n,d}$ converges to 0 as $m, n\rightarrow\infty$. Thus we have consistency under both $H_0$ and $H_1$. \qed
\end{enumerate}

\bibliographystyle{Chicago}
\bibliography{BERET_BIB.bib}

\end{document}